\documentclass[12pt]{amsart}
\usepackage{amscd,amssymb,graphics}
\usepackage{mathrsfs} 
\newtheorem{theorem}{Theorem}[section]
\newtheorem{corollary}[theorem]{Corollary} 
\newtheorem{lemma}[theorem]{Lemma}
\newtheorem{proposition}[theorem]{Proposition}
\theoremstyle{definition}
\newtheorem{definition}[theorem]{Definition}

\theoremstyle{remark}
\newtheorem{remark}[theorem]{Remark}

\newtheorem{example}[theorem]{Example}

\numberwithin{equation}{section}

\newcommand{\abs}[1]{\lvert#1\rvert}
\newcommand{\vNa}{von Neumann algebra}
\def\norm#1{\left\Vert#1\right\Vert}
\def\R {{\mathbb R}}
\def\Q {{\mathbb Q}}
\def\I {{\mathbb I}}
\def\C {{\mathbb C}}
\def\N{{\mathbb N}}
\def\e{{\varepsilon}}
\def\Z {{\mathbb Z}}
\def\Un{{\mathcal U}}
\def\H{{\mathscr H}}

\def\B{{\mathscr B}}

\def\supp{{\mathrm{supp}}}
\def\Ad{{\mathrm{Ad}}}
\def\s{{\mathbb S}}
\def\1{{\mathbf 1}}
\def\UCB{{\mathrm{UCB}}\,}
\def\Aut{{\mathrm{Aut}}}
\def\Cnt{{\mathrm{Cnt}}}

\def\Homeo{{\mathrm{Homeo}}}
\def\Inn{{\mathrm{Inn}}}
\def\Tr{{\mathrm{Tr}}}
\def\Ad{{\mathrm{Ad}\,}}
\def\Iso{{\mathrm{Iso}}}
\def\RUCB{{\mathrm{RUCB}}\,}
\def\LUCB{{\mathrm{LUCB}}\,}
\def\Bij{{\mathrm{Bij}}}

\def\HL{{\underline{H}_1{\mathcal L}\iota}}

\begin{document}

\title[Some extremely amenable groups]{Some extremely amenable groups
related to operator algebras and ergodic theory}

\author[T. Giordano]{Thierry Giordano}
\author[V. Pestov]{Vladimir Pestov}
\address{Department of Mathematics and Statistics, 
University of Ottawa, 585 King Edward Ave., Ottawa, Ontario, Canada K1N 6N5.
}
\email{giordano@uottawa.ca, vpest283@uottawa.ca}




\begin{abstract} 
A topological group $G$ is called extremely amenable
if every continuous action of $G$ on a compact space has a fixed point.
This concept is linked with geometry of high dimensions
(concentration of measure). 
We show that a
von Neumann algebra is approximately finite-dimensional
if and only if its unitary group with the strong
topology is the product of an extremely amenable group with a compact group,
which strengthens a result by de la Harpe.
As a consequence, a $C^\ast$-algebra $A$ is nuclear if and only if
the unitary group $U(A)$ with the relative weak topology is strongly
amenable in the sense of Glasner. We prove that
the group of automorphisms of a Lebesgue space with a non-atomic
measure is extremely amenable
with the weak topology and establish
a similar result for groups of non-singular
transformations. As a consequence,
we prove extreme amenability of the groups of isometries of $L^p(0,1)$,
$1\leq p<\infty$, extending a classical result of Gromov and Milman ($p=2$).
We show that a measure class
preserving equivalence relation $\mathcal R$
on a standard Borel space is amenable
if and only if the full group $[{\mathcal R}]$, equipped with the uniform
topology, is extremely amenable. 
Finally, we give natural examples of concentration to a nontrivial space
in the sense of Gromov occuring in the automorphism groups of injective factors
of type $III$.
 \end{abstract}

\maketitle


\section{Introduction}
Extreme amenability \cite{Gra}
(or the fixed point on compacta property \cite{Mit})
is a relatively recent concept, useful for study of `large' 
(non locally compact) topological groups.
Recall that a topological group $G$ is {\it amenable} \cite{dlH} if every
continuous action of $G$ by affine transformations on a convex compact subset
of a topological vector space has a fixed point. 
A topological group $G$ is {\it extremely
amenable} if every continuous action of $G$ on a compact space 
has a fixed point.

This is a strong nonlinear fixed point property,
never observed in locally compact groups \cite{GraL,V}.
Beginning with the discovery by Gromov and Milman \cite{GrM} that 
the unitary group of an infinite
dimensional Hilbert space with the strong operator topology is
extremely amenable, a large number of 
concrete large groups of importance are now known to have this property.

The concept is interesting for several reasons. 
Extreme amenability 
is closely linked to asymptotic geometric analysis,
and in many examples it is a manifestation of
the phenomenon of concentration of measure on high-dimensional structures,
captured in the concept of a L\'evy group
\cite{GrM,Gl,GP}.
Extreme amenability of groups of automorphisms of various
combinatorial structures is intimately related
to Ramsey-type theorems for those structures \cite{P1,KPT}. 
While surely not all
infinite dimensional groups are extremely amenable, the fixed point on
compacta property
of subgroups can be used to make conclusions about the dynamical
properties of larger groups, for instance, to describe their universal minimal
compact flows \cite{P1,Gl-W1}. 
And the recent work by Glasner, Tsirelson and Weiss \cite{GTW} and
Glasner and Weiss \cite{Gl-W3}
links L\'evy groups to the problem of (non)existence of spatial models for
near actions.

No non-trivial locally compact group has the fixed point on
compacta property, which is thus exclusively
a property of `large' topological groups. This was shown by Granirer and Lau
\cite{GraL}. Later Veech \cite{V} has proved his well-known theorem
stating that every locally compact group acts freely on a
compact space. 
(Now there are at least four known proofs of this result, see 
\cite{KPT} and references therein.)

The first example of a non-trivial
extremely amenable group has been constructed by Herer and Christensen
\cite{HC}.
Gromov and Milman \cite{GrM} have shown that the unitary group of an infinite
dimensional Hilbert space, equipped with the strong operator topology, is extremely amenable.
The idea of the proof, using concentration of measure on
high-dimensional structures, has turned out to be applicable in a great
many situations. For example, Glasner, and also (independently, unpublished)
Furstenberg and Weiss, have shown that the group $L^0(\I,U(1))$
of measurable maps from
the standard Lebesgue space to the circle rotation group $U(1)$, equipped with
the topology of convergence in measure, is extremely amenable, see
\cite{Gl}. The group $U(1)$ here
can be replaced by any amenable locally compact group \cite{P5}.

A direct link between extreme amenability and 
Ramsey theory has been established 
by one of the present authors, who has proved in \cite{P1}
extreme amenability of the group $\Aut(\Q,\leq)$ of order-preserving
bijections of the rationals, equipped with the topology of pointwise
convergence on $\Q$ viewed as discrete. This statement is a reformulation
of the classical Finite Ramsey Theorem. The recent paper by Kechris,
Pestov and Todorcevic \cite{KPT} explores this trend extensively, by establishing
extreme amenability of groups of automorphisms of numerous countable
Fra\"\i ss\'e structures.
It should be noted that extreme amenability of
a topological group $G$ can be restated in terms of the so-called
Ramsey-Dvoretzky-Milman property (Section 9.3 in \cite{Gr1})
of transitive isometric actions of $G$ on
metric spaces, which is in turn linked to both discrete Ramsey theory and 
Ramsey-type properties in geometric functional analysis \cite{P6}.

The universal minimal flow $M(G)$ of a topological group $G$ is, 
in the case where $G$ is locally compact and non-compact, a
highly non-constructive object (for instance, 
non-metrizable, \cite{KPT}, Appendix 2). 
Of course, if $G$ is an extremely amenable group, then
$M(G) =\{\ast\}$ is a singleton. 
It appears that the first instance where the flow $M(G)$,
different from a point,
was described explicitely, was the case where $G={\mathrm{Homeo}}_+(\s^1)$ is the
group of orientation-preserving homeomorphisms of the circle with the
compact-open topology. In this situation, $M(G)$ is the circle
$\s^1$ itself, equipped with the canonical action of $G$ 
\cite{P1}. The proof was using the existence of a large
extremely amenable subgroup, namely ${\mathrm{Homeo}}_+(\I)$. 
Glasner and Weiss 
have subsequently described in \cite{Gl-W1} the universal minimal flow
of the infinite symmetric group with its (unique) Polish topology as the
space of linear orders on the natural numbers, and
then in \cite{Gl-W2} the universal minimal
flow of the group of homeomorphisms of the Cantor set, $C$, as the space of
maximal chains of closed subsets of $C$ (a construction proposed by
Uspenskij \cite{Usp00}). 
Numerous other examples of explicit computations of universal
minimal flows of groups of automorphisms of countable structures can be
found in the above cited paper \cite{KPT}.

Overall, the papers \cite{P6} and \cite{KPT} will together provide an
introduction to, a survey of, and
bibliographical references to most 
of what is known to date about extremely amenable groups. 
An account of basic ideas of the theory is provided in the
set of lecture notes by one of the present authors \cite{P05}.

In this article we will 
establish extreme amenability of some concrete topological groups of 
importance using the concentration of measure, as
well as connect the property with known concepts from
operator algebras and ergodic theory.
We investigate extreme amenability of topological groups
of two types: unitary groups of von Neumann algebras and group of
transformations of spaces with measure. 

It was proved by de la Harpe \cite{dlH1} that a von Neumann algebra $M$ is
approximately finite dimensional if and only if its unitary group $U(M)$,
equipped with the relative $s(M,M_\ast)$-topology, is amenable. 
We strengthen this result by showing that this is the case if and only if
$U(M)$ is the product of a compact group and an extremely amenable group.
Paterson \cite{Pat}
has deduced from de la Harpe's result the following: a $C^\ast$-algebra
$A$ is nuclear if and only if the unitary group $U(A)$, equipped with the
relative weak (that is, $s(A,A^\ast)$) topology, is amenable.
We describe nuclear $C^\ast$-algebras as those $C^\ast$-algebras $A$ 
whose unitary group $U(A)$ with the relative weak topology is strongly
amenable in the sense of Glasner \cite{Gl1}. 

Let now $(X,\nu)$ be a standard non-atomic Borel probability space. We show
that the group $\Aut(X,\nu)$ of all measure-preserving automorphisms
of $(X,\nu)$, equipped with the weak topology, is extremely amenable.
This result is a consequence of the Rokhlin--Kakutani lemma and the
concentration of measure on finite symmetric groups discovered by Maurey
\cite{Ma}.
The group $\Aut(X,\nu)$, equipped with the uniform topology, is no longer
extremely amenable. We also establish the extreme amenability of the
group $\Aut^\ast(X,\nu)$ of measure class preserving transformations of
$(X,\nu)$, equipped with the weak topology. As a consequence of this result
and a description of groups of isometries of spaces $L^p(X,\mu)$ belonging to
Banach \cite{banach} and Lamperti \cite{lamperti}, we
prove that those groups, equipped with the strong operator topology, 
are extremely amenable for all $1\leq p<\infty$, 
extending Gromov and Milman's classical result \cite{GrM} (corresponding to $p=2$).

We consider measure class preserving equivalence relations $\mathcal R$ on the
standard Lebesgue measure space and prove that such a relation is amenable
in the sense of Zimmer \cite{zimmer-inv}, \cite{zimmer-jfa} 
if and only if the full group $[{\mathcal{R}}]$, equipped with the uniform
topology, is extremely amenable. In order to obtain this and the previous
result, we generalize Maurey's theorem to automorphism groups of 
measure spaces with finitely many points.

The concentration of measure phenomenon can be interpreted in terms of
convergence of a family of spaces with metric and measure ($mm$-spaces)
to a one-point space with respect to a suitable metric, as was shown by
Gromov \cite{Gr}. This leads to a more general
concept of concentration to a non-trivial
space. In the concluding part of the article, we show some natural
examples of concentration of compact subgroups and other subobjects of the
groups of automorphisms of injective factors of type $III$ to non-trivial
spaces.

Some of the above results have been announed by the present authors in
\cite{GP}.

\subsection*{Acknowledgements}
The investigation has been partially supported by NSERC
operating grants and by the Swiss Fund for Scientific Research (both authors), 
as well as the Marsden Fund of the Royal Society of
New Zealand and a University of Ottawa start-up grant (V.P. only).
Both authors are grateful to Alexander Kechris for useful comments, and
to Pierre de la Harpe for his hospitality at the University of Geneva where the
present research collaboration started in April 1999. 
The authors are also grateful to the referee  whose suggestions have led to a
considerable reworking of the article,
whose several proofs have been simplified and made clearer.

\section{Concentration and extreme amenability}

In this section we will outline
a general scheme of deducing fixed point theorems from the
phenomenon of  concentration of measure,
introduced by Gromov and Milman \cite{GrM}. We work in a slightly more
general context and prove some new results, in particular a concentration of
measure result for automorphism groups of finite measure spaces. 

\subsection{Concentration of measure and L\'evy families}

For a subset $A\subseteq X$ of a separated uniform space $X=(X,{\mathcal U})$
let $V[A]=\{x\in X \colon \exists a\in A, ~(x,a)\in V\}$ be the
$V$-neighbourhood of $A$.
We say that a net $(\mu_\alpha)$ of probability measures on $X$
has the {\it L\'evy concentration property,} or simply
{\it concentrates}, if for every family of Borel subsets
$A_\alpha\subseteq X$ satisfying
$\liminf_\alpha\mu_\alpha(A_\alpha)>0$
and every entourage $V\in{\mathcal{U}}_X$ one has
$\mu_\alpha(V[A_\alpha])\to_{\alpha} 1$.

A triple $(X,d,\mu)$, where $d$ is a metric on a set $X$ and $\mu$
is a probability Borel measure on the metric space $(X,d)$, is called a
{\it metric space with measure}, or else an $mm${\it -space}
\cite{Gr}. An infinite family $(X_n,d_n,\mu_n)$ of $mm$-spaces is a
{\it L\'evy family} if, whenever $A_n\subseteq X_n$ are Borel subsets satisfying
$\liminf_n\mu_n(A_n)>0$, one has for each $\e>0$ $\lim_{n\to\infty}
\mu_n((A_n)_\e)=1$, where $A_\e$ denotes the $\e$-neighbourhood of a set $A$. 

The above property is equivalent to
concentration of the family $(\mu_n)$ considered
as probability measures supported on the disjoint union 
$\oplus_{n=1}^\infty X_n$, equipped
with a metric $d$ inducing the metrics $d_n$ on each $X_n$ and making $X_n$ into
an open and closed subset. 

The next Lemma first appears in \cite{GrM}, 2.2. For a proof, see
Proposition 3.8 in \cite{L}. (The case of uniform spaces clearly
reduces to that of $mm$-spaces.) 

\begin{lemma}
Let a net of probability measures $(\mu_\alpha)$
(resp. $(\nu_\alpha)$) on an uniform space $X$ (resp. $Y$) concentrate.
Then the net of product measures $(\mu_\alpha\otimes \nu_\alpha)$
concentrates on the product space $X\times Y$.
\label{subadd}
\qed
\end{lemma}


The {\it concentration function}, $\alpha_X$, of an $mm$-space \cite{GrM} is defined
for $\e\geq 0$ by
\[
\alpha_X(\e)=\begin{cases}\frac 12, & \text{if $\e=0$,} \\
1-\inf\{\mu(A_\e)\colon \text{$A\subseteq X$ is Borel, $\mu(A)\geq\frac 12$}\},
&\text{if $\e>0$.}
\end{cases}
\]

\begin{definition} A family $(X_n,d_n,\mu_n)$ of $mm$-spaces is a
{\em normal L\'evy family} if for some constants $C_1,C_2>0$ 
\[\alpha_{X_n}(\e)\leq C_1 \exp(-C_2 n \e^2).\]
\end{definition}

The following is simple and well-known.

\begin{lemma}
\label{above}
Let $(Y,\nu)$ be a probability measure space, and let $d,\rho$ be two
measurable metrics on $Y$ such that the identity map
\[(Y,d,\nu)\to (Y,\rho,\nu)\]
is Lipschitz with constant $L$. 
Denote by $\alpha_d$ the concentration function
of the $mm$-space $(Y,d,\nu)$, and similarly by $\alpha_\rho$ 
the concentration function of the $mm$-space $(Y,\rho,\nu)$.
Then for every $\e>0$
\[\alpha_\rho(\e)\leq \alpha_d(L^{-1}\e).\]
\end{lemma}

\begin{proof}
Indeed, the $L^{-1}\e$-neighbourhood (w.r.t. $d$) of an arbitrary
measurable subset $A\subseteq Y$ having measure
$\geq \frac 12$ is contained in the $\e$-neighbourhood (w.r.t. $\rho$)
of $A$, and therefore the measure of the latter is at least as large as
the measure of the former. 
\end{proof}

A family $(X_n,d_n,\mu_n)$ is a
L\'evy family if and only if the concentration functions
$\alpha_n$ of $X_n$ converge to zero pointwise as $n\to\infty$ on the
interval $(0,+\infty)$.

A well-known manifestation of concentration of measure
(see \cite{M88,MS,L}) is that on a highly-concentrated space every
uniformly continuous function is ``nearly constant nearly everywhere.''
More precisely:

\begin{lemma} 
\label{nearlyconstant}
Let a family $(\mu_\beta)$ of probability measures on a
uniform space $X$ be L\'evy. Then for every bounded uniformly
continuous function $f$ on $X$ there exists a net of constants
$(c_\beta)$ such that for every $\e>0$,
\begin{equation}
\mu_\beta\{x\in X\colon \abs{f(x)-c_\beta}>\e\}\to_{\beta} 0.
\label{nearly}
\end{equation}
\end{lemma} 

\begin{proof}[Proof, sketch] One can choose as 
$c_\beta$ median values for $f$ with respect to $\mu_\beta$, 
that is, numbers with the property
that sets $M_{\beta}^-(f)=\{x\in X\colon f(x)\leq c_\beta\}$ and
$M_{\beta}^+(f)=\{x\in X\colon f(x)\geq c_\beta\}$ have 
both $\mu_\beta$-measure at least $1/2$. Let $V=V(\e)\in\Un_X$
be such that $\abs{f(x)-f(y)}<\e$ whenever $(x,y)\in V$.
The set in Eq. (\ref{nearly}) is
contained in the intersection of complements to $V$-neighbourhoods of
$M^{\pm}_\beta$. 
Now the upper bound is obtained by
applying the definition of a L\'evy family of measures.
\end{proof}

Here are a few examples of L\'evy families.

\begin{example} The family of special unitary groups $SU(n)$, $n\in\N$,
equipped with the normalized Haar measure and the Hilbert-Schmidt metric
is a L\'evy family \cite{GrM}.
\end{example}

\begin{example} 
\label{prod}
Let $X=(X,\mu)$ be a probability measure space. Introduce on $X^n$
the product measure $\mu^{\otimes n}$ and the normalized Hamming distance
$d(x,y) =\frac 1n\vert\{i\colon x_i\neq y_i\}\vert$. The family $X^n$, $n\in\N$, is then
a L\'evy family, with the concentration function $\alpha_n$ of $X^n$
satisfying
\begin{equation}
\alpha_n(\e)\leq 2\exp(-\e^2n)
\label{concfuncprod}
\end{equation}
(See e.g. \cite{Ta}, Prop. 2.1.1.)

Here is an immediate consequence. Let $X_1,X_2,\ldots$ be a sequence of
probability measure spaces. Equip the $n$-fold product $Y_n = X_1\times\ldots
\times X_n$ with the product measure and the normalized Hamming distance. Then the
family $(Y_n)_{n=1}^\infty$ is L\'evy, with the same concentration functions
as in Eq. (\ref{concfuncprod}). (Apply Example \ref{prod} and
Lemma \ref{above} to
the L\'evy family $\I^n$ with Lebesgue measures and the normalized Hamming distance.)
\end{example}

\begin{example} Symmetric groups ${\mathfrak{S}}_n$ of rank $n$, equipped
with the {\em normalized Hamming distance}
\[d_{n}(\sigma,\tau):=\frac 1n \abs{\{i\colon \sigma(i)\neq\tau(i)}\]
and the normalized counting measure
\[\mu(A):=\frac{\abs{A}}{n!}\]
form a normal L\'evy family, with the concentration functions
satisfying the estimate
\begin{equation}
\alpha_{{\mathfrak{S}}_n}(\e)\leq \exp(-\e^2n/32).
\label{maureyest}
\end{equation}
(Maurey \cite{Ma}; see also \cite{MS}, Theorem 7.5, and
\cite{L}, Corollary 4.3.)
\label{maurey}
\end{example}

There are many known techniques for establishing concentration inequalities. 
The following approach, based on martingales, allows one 
to prove Examples \ref{prod} and \ref{maurey}, as well as to obtain some useful
variations and generalizations.

Following Milman and Schechtman \cite{MS}, say that a finite metric space
$(\Omega,d)$ is {\it of length at most} $\ell$ if there exist positive numbers
$a_1,a_2,\ldots,a_n$ with $\ell = (\sum_{i=1}^n a_i^2)^{1/2}$ and a refining
sequence
$(\Omega^k)_{k=0}^n$, $\Omega^k = \{A_i^k\}_{i=1}^{m_k}$ of partitions of 
$\Omega$ such that $\Omega^0=\{\Omega\}$, 
$\Omega^n = \{\{x\}\colon x\in\Omega\}$,
and for every $k=1,\ldots,n$ and every two elements $A,B$ of the partition
$\Omega^k$ that are both contained in some element of $\Omega^{k-1}$
there exists a bijection $\varphi\colon A\to B$ with
$d(x,\varphi(x))\leq a_k$ for all $x\in A$. 

The following is Theorem 7.8 in \cite{MS}, and the present improved
constants are taken from Theorem 4.2 in \cite{L}.

\begin{theorem} Let $(\Omega,d)$ be a finite metric space of length
at most $\ell$, equipped with the normalized counting measure.
Then
\[\alpha_{\Omega}(\e)\leq \exp\left(- \frac{\e^2}{8\ell^2} \right).\]
\qed
\label{length}
\end{theorem}

The following is Theorem 7.12 in \cite{MS}, and (with the present, better
constants) Theorem 4.4 in \cite{L}.

\begin{theorem} Let $G$ be a compact group equipped with a bi-invariant
metric $d$, and let $\{e\}=H_0<H_1<\ldots<H_n=G$ be a sequence of closed
subgroups. 
Equip every factor-space $H_{i+1}/H_i$ with the factor-distance of $d$, and
let $d_i$ denote the diameter of $H_{i+1}/H_i$. Then the concentration
function of the $mm$-space $(G,d,\mu)$, where $\mu$ is the normalized Haar
measure, satisfies
\[\alpha_G(\e)\leq \exp\left(- \frac{\e^2}{8\sum_{i=0}^{n-1} {d_i}^2} 
\right).\]
\label{forgroups}
\qed
\end{theorem}

\subsection{Two new examples of L\'evy families} Here we apply
Theorems \ref{length} and \ref{forgroups} to deduce two new natural examples of
L\'evy families, which we will use in the sequel. 


\begin{corollary} Let $(X,\mu)$ be a probability measure space with
finitely many points, and let $(K,d)$ be a compact metric group of
diameter $1$. The concentration function of the compact metric group
$L^1(X,\mu;K)$, equipped with the normalized Haar measure and the 
$L^1(\mu)$-metric 
\[d^\dag(f,g) = \int_{X} d(f(x),g(x))~d\mu(x),\] 
satisfies
\[\alpha(\e)\leq \exp\left(- \frac{\e^2}{8\sum_{x\in X} \mu(\{x\})^2} 
\right).\]
\label{weneed}
\end{corollary}

\begin{proof} If $X=\{x_1,x_2,\ldots,x_n\}$, then let $H_i$ consist of all
functions vanishing on $x_{i+1},x_{i+2}$, $\ldots,x_n$. The group
$H_{i+1}/H_i$ is isometrically isomorphic to the group $(K,\mu(\{x_{i+1}\})d)$,
whose diameter equals $\mu(\{x_{i+1}\})$, and the result follows by
Theorem \ref{forgroups}.
\end{proof}

\begin{corollary} Let $(X_n,\mu_n)$ be a sequence of probability
spaces with finitely many points in each, having
the property that the mass of the largest atom in $X_n$ goes
to zero as $n\to\infty$. Let $(K,d)$ be a compact metric group of diameter $1$.
Then the family of compact groups $L^1(X_n,\mu_n; K)$, equipped with the
normalized Haar measures and the 
$L^1(\mu)$-metrics, is L\'evy, with
\[\alpha_{X_n}(\e)\leq  \exp
\left(-\frac{\e^2}{8\max_{x\in X_n}\mu_n(\{x\})}\right).\]
\label{largestatom1}
\end{corollary}

\begin{proof} Let $v_n=\left(\mu_n(\{x\})\right)_{x\in X_n}$.
By H\"older's inequality,
\[ \sum_{x\in X_n} \mu_n(\{x\})^2 = \langle v_n,v_n\rangle \leq
\norm{v_n}_1\cdot \norm{v_n}_\infty = \max_{x\in X_n}\mu_n(\{x\})
\to 0\mbox{ as } n\to\infty.\]
Corollary \ref{weneed} now applies.
\end{proof}

Let $X= (X,\mu)$ be a measure space. Denote by $\Aut^\ast(X,\mu)$ the group of
all invertible measurable and
non-singular transformations of $X$. In particular, if $X$ has finitely
many points and the measure has full support, then
$\Aut^\ast(X,\mu)$ is, as an abstract group, just the symmetric group of rank
$n=\abs X$. 

Define a left-invariant metric on 
$\Aut^\ast (X,\mu)$ (the uniform metric) as follows:
\begin{equation}
\label{uniformmetric}
d_{unif}(\tau,\sigma)=\mu\{x\in X\colon \tau(x)\neq\sigma(x)\}.
\end{equation}
For instance, if a finite set $X$ is equipped with a uniform measure,
the uniform metric on $\Aut^\ast(X,\mu)\cong {\mathfrak{S}}_n$ coincides with the Hamming
distance. For every measure space, the uniform metric generates a group
topology, known as the {\it uniform topology}, which will play an important
role later (cf. Subsection \ref{twotopologies}).

\begin{theorem} Let $X=(X,\mu)$ be a probability space with finitely many
points. The concentration function $\alpha$ of the
automorphism group $\Aut^\ast(X,\mu)$, equipped with the uniform metric
(\ref{uniformmetric}) and the normalized counting measure, satisfies
\begin{equation}
\alpha(\varepsilon) \leq  
\exp\left(-\frac{\e^2}{32\sum_{x\in X}\mu(\{x\})^2}\right).
\end{equation}
\label{concauto}
\end{theorem}

\begin{proof} 
Let $X = \{x_1,x_2,\ldots,x_n\}$, where
\[\mu(\{x_1\}) \geq \mu(\{x_{2}\}) \geq \ldots\geq \mu(\{x_{n-1}\})
\geq \mu(\{x_{n}\}). \]
For $k=0,1,\ldots,n$, let $H_k$ be the subgroup stabilizing each element 
$x_1,\ldots,x_k$. 
Thus,
\[\Aut^\ast(X,\mu) = H_0 >H_1>H_2>\ldots>H_n = \{e\}.\]
Let $\Omega^k$ be the partition of $G = \Aut^\ast(X,\mu)$ into left
$H_k$-cosets $\sigma H_k$, $\sigma\in G$. 

Suppose $A=\sigma H_k$ and $B=\tau H_k$ are contained in the same left
$H_{k-1}$-coset. Then $\sigma(x_i) =\tau(x_i)=a_i$ for $i=1,2,\ldots,k-1$,
while $a=\sigma(x_k)$ and $b=\tau(x_k)$ need not coincide. 
Thus, elements $\pi$ of $A$ are defined by the conditions
\[\pi(x_1) = a_1,~ \pi(x_2)=a_2,~\ldots,~\pi(x_{k-1})=a_{k-1},~ \pi(x_k)=a,\]
while
elements $\pi\in B$ are defined by the conditions
\[\pi(x_1) = a_1,~ \pi(x_2)=a_2,~\ldots,~\pi(x_{k-1})=a_{k-1},~ \pi(x_k)=b.\]

Let $t_{a,b}$ denote the transposition of $a$ and $b$. 
Consider the map 
\begin{equation}
\label{ourmap}
\varphi\colon\sigma H_k\ni j\mapsto t_{a,b}\circ j \in \tau H_k.
\end{equation}
Clearly, $\varphi$ is a bijection between $A=\sigma H_k$ and $B=\tau H_k$.
The values of
$j$ and $t_{a,b}\circ j$ differ in at most two inputs, 
$x_k=j^{-1}(a)$ and $j^{-1}(b)$. 
Since $a,b\notin \{a_1,a_2,\ldots,a_{k-1}\}$, it follows that $j^{-1}(a),
j^{-1}(b)\notin
\{x_1,x_2,\ldots,x_{k-1}\}$ and 
$\mu(\{j^{-1}(b)\})\leq \mu\{(x_k\})$.
We conclude: for every $j\in A=\sigma H_k$
\[d_{unif}(j,\varphi(j))\leq 2 \mu(\{x_k\}).\]
Consequently, the metric space $\Aut^\ast(X,\mu)$ has length at most
$\ell=2\left(\sum_{i=1}^n \mu(\{x_i\}^2)\right)^{1/2}$, and 
Theorem \ref{length} accomplishes the proof.
\end{proof}

\begin{corollary} Let $(X_n,\mu_n)$ be a sequence of probability
spaces with finitely many points in each, having
the property that the mass of the largest atom in $X_n$ goes
to zero as $n\to\infty$. Then the family of automorphism groups
$\Aut^\ast(X_n,\mu_n)$, $n\in\N$, equipped with the uniform metric and the
normalized counting measure, is L\'evy, with
\[\alpha_{X_n}(\e)\leq  \exp
\left(-\frac{\e^2}{32\max_{x\in X_n}\mu_n(\{x\})}\right).\]
\label{largestatom}
\end{corollary}

\begin{proof} Same as the proof of Corollary \ref{largestatom1}.
\end{proof}

\subsection{Fixed points}
If $X=(X,{\mathcal U})$ is a uniform space, then $\sigma X$ will denote
the {\it Samuel}, or {\it universal, compactification} of $X$, that is, 
the space of maximal ideals of the
$C^\ast$-algebra $\UCB(X)$ of all uniformly continuous bounded complex-valued
functions on $X$. (See \cite{E}.)

If a group $G$ acts on $X$ by uniform isomorphisms, this action uniquely
extends to an action of $G$ on $\sigma X$ by homeomorphisms. 

The following result is essentially Th. 2.12 in \cite{P5}, to which we  
give a sketch of a different proof.

\begin{theorem}
 Let a group $G$ act by uniform isomorphisms on a uniform
space $X$. Suppose a net $(\nu_\alpha)$ of probability measures on $X$
concentrates and has the property:
$(\ast)$ for every $g\in G$, every $A\subseteq X$ and
every $V\in{\mathcal{U}}_X$, 
\[\limsup_\alpha \left\lvert \nu_\alpha(V[A])-
\nu_\alpha(g\cdot V[A])\right\rvert<1.\]
Then the Samuel compactification $\sigma X$ has a $G$-fixed point.
\label{two-twelve}
\end{theorem}

\begin{proof}[Proof, sketch] The net $(\phi_\alpha)$ of means on $C(\sigma X)$, defined by
$\phi_\alpha(f) = \int_X f(x) d\nu_\alpha(x)$, has a weak$^\ast$ cluster point,
say $\phi$. By selecting a subnet, assume that $\phi_\alpha\to\phi$
in weak$^\ast$ topology. As a consequence of Lemma \ref{nearlyconstant}, 
for any bounded uniformly continuous function $f$ on $X$, the 
$L_1(\nu_\alpha)$-distance between $f$ and a suitable constant function
(for instance, $M_\alpha(f)$, a $\nu_\alpha$-median of $f$) goes to
zero as $\alpha\to \infty$. Consequently, the $L_1(\nu_\alpha)$-distance
between $f\cdot g$ and $M_\alpha(f)\cdot M_\alpha(g)$ tends to zero as well,
and
\[\lvert{\int fg d\nu_\alpha - \left(\int f d\nu_\alpha\right)\cdot
\left(\int g d\nu_\alpha\right)}\rvert\to_\alpha 0, \]
for all $f,g\in C(\sigma X)$, implying that the mean $\phi$ is multiplicative.

The condition $(\ast)$, jointly with concentration of $(\nu_\alpha)$,
easily implies that whenever $\liminf_\alpha \nu_\alpha (A)>0$,
$g\in G$, and $V\in{\mathcal{U}}_X$, one has
$\limsup_\alpha \nu_\alpha (VA\cap gVA)=1$. Consequently, if $f\in C(\sigma X)$, 
there is a subnet $(\beta)$ with $M_\beta(f)-M_{\beta}(^gf)\to 0$, where
$^gf(x)=f(g^{-1}x)$. We
conclude: $\phi_{\beta}(f) -\phi_{\beta}(^gf) \to_\beta 0$, and as
$\phi_{\beta}(f)\to_\beta\phi(f)$, it follows that $\phi$ is a $G$-invariant
mean.
\end{proof}

Let $X$ be a topological group $G$ equipped with the
right uniform structure, which has as a basis the entourages of the diagonal 
of the form $V_R=\{(x,y)\in G\times G \mid xy^{-1}\in V\}$,
where $V$ is a neighbourhood of $e$ in $G$. Notice that
for a subset $A$ of $G$, the uniform neighbourhood $V_R[A]$ is just the
product $VA$.

The corresponding Samuel compactification,
denoted by $S(G)$, is the {\it greatest ambit} of
$G$, that is, the maximal compact $G$-space containing a distinguished point with 
an everywhere dense orbit. 
The existence of a $G$-fixed point in $S(G)$ is therefore equivalent to
extreme amenability of $G$.

The condition ($\ast$) is a weak form of invariance
of the measures $\nu_\alpha$. In many
examples, one can assume that the measures $\nu_\alpha$ are 
virtually invariant, that is, for every
$g\in G$, one has $g\ast \nu_\alpha = \nu_\alpha$ for all $\alpha \geq \alpha_0$. 
A topological group $G$ is called a {\it L\'evy group} \cite{GrM,Gl}
if it contains a family $\mathcal K$ of compact 
subgroups, directed by inclusion, having everywhere dense union in $G$, and
such that the corresponding normalized 
Haar measures, $\mu_K$, on the groups 
$K\in{\mathcal{K}}$ concentrate in $G_R$. 

The presently known examples of L\'evy groups are numerous. The one of the greatest importance for us is the original example essentially obtained in \cite{GrM}.

\begin{example}[Gromov and Milman]
The group of all unitary operators of the form $\I+C$, where $C$ is a compact operator
of Schatten class $2$, equipped with the Hilbert-Schmidt metric, is a L\'evy group.

It follows at once from the fact that the
increasing chain of subgroups $SU(n)$, $n\in\N$ forms a L\'evy family with regard to the Hilbert-Schmidt metric (\cite{GrM}, Example 3.4), and its union is well-known to be dense in the operator group in question.
\label{ex:ugr}
\end{example}

We expand the concept of a L\'evy group as follows.

\begin{definition}
\label{glg}
Let us say that a topological group $G$ is a 
{\it generalized L\'evy group} if there is a net 
$(K_\alpha)_{\alpha\in A}$
of compact subgroups of $G$ with the following properties:
\begin{enumerate}
\item The family of normalized Haar measures $\mu_\alpha$ on
$K_\alpha$ concentrates in $G$;
\item For every non-empty open subset $V\subseteq G$ there is an $\alpha\in A$
such that $V\cap K_\beta \neq\emptyset$ for all $\beta\geq\alpha$.
\end{enumerate}
\end{definition}

\begin{remark} The second condition is equivalent to the
following, formally stronger, condition: 
for every finite collection $g_1,g_2,\dots,g_N\in G$, $N\in\N$, and for 
every neighbourhood $V$ of the identity in $G$ 
there is an $\alpha\in A$ such that for all
$\beta\geq\alpha$ and $i=1,2,\dots,N$, one has $Vg_i\cap K_\beta \neq\emptyset$.
\end{remark}

\begin{theorem} Every generalized L\'evy group is extremely amenable.
\qed
\label{gm}
\end{theorem}

\begin{proof} In view of Theorem \ref{two-twelve}, it is enough to
verify the condition ($\ast$). Assume it fails. Then for some $g\in G$, some
(non-empty) $A\subset G$, and some symmetric neighbourhood
of identity, $V$, one can construct a subnet $(K_{\beta})$ of
groups with the properties
$\mu_\beta(V^2A)\to_\beta 1$ and $\mu_\beta(gV^2A)\to_\beta 0$.

Because of condition (2) in the definition of a generalized L\'evy group,
one can assume without loss in generality that for all $\beta$,
$gV\cap K_\beta\neq\emptyset$. Choose a net
$g_\beta'\in K_\beta$
with $g=g_\beta'v$, where $v\in V$. Since $vV^2A\supseteq VA$, one has
$\nu_\beta(gV^2A)\geq \nu_\beta(g_\beta'VA) = \nu_\beta(VA)\to_\beta 1$, meaning that
$\nu_\beta(gV^2A)\to_\beta 1$ as well, a contradiction.
\end{proof}

\begin{corollary}[Gromov and Milman]
\label{everylevy}
Every L\'evy group is extremely amenable.
\qed
\end{corollary}

To deduce Corollary \ref{everylevy} directly 
from Theorem \ref{two-twelve}, it is enough to note the extreme 
amenability of the
everywhere dense subgroup $\cup{\mathcal{K}}$ of $G$, since both topological
groups share the same greatest ambit.

There are still situations where one needs the condition 
($\ast$) in its full generality. 
For instance, this is how one proves the following Theorem \ref{gener}.

Let $G$ be a Hausdorff topological group and let $(X,\mu)$ be a Lebesgue space. 
Denote by $L(X,\mu;G)$, or simply $L(X,G)$
the group of equivalence classes modulo $\mu$ of all (strongly) measurable maps
(cf. \cite{edwards}, 8.14.1(b))
from $(X,\mu)$ to $G$, equipped with pointwise multiplication. 
The {\em topology of convergence in measure} on $L(X,G)$ has, as basic
neighbourhoods of identity sets of the form
\[
[V,\e] :=\{ g\in L(X,G) \colon \mu\{x\in X\colon g(x)\notin V\}<\e\},
\]
where $V$ is a neighbourhood of $e_G$ in $G$ and $\e>0$.
The topology of convergence in measure is a Hausdorff group topology
on $L(X,G)$,
and we will denote the resulting topological group by $L^0(X,\mu;G)$
or $L(X,G)$.

\begin{theorem}
Let $G$ be an amenable locally compact group and let $X$ be a
non-atomic Lebesgue measure space. 
Then the group $L^0(X,G)$ of all measurable maps from $X$ to $G$, 
equipped with the topology of convergence in measure, is extremely amenable. 
\label{gener}
\end{theorem}

\begin{proof}[Proof, sketch]
By force of amenability of $G$, there is an invariant probability measure,
$\nu$, on the greatest ambit $S(G)$ of $G$. Choose 
a net $(\nu_\alpha)$ of probability measures on $G$ (for instance, with
finite support), weak$^\ast$ converging to $\nu$. For every measurable
partition $\mathcal W$ of $X$ define a measure $\nu_{\alpha,{\mathcal W}}$
on $L^0(X,G)$, supported on the subgroup of all $\mathcal W$-simple maps
$X\to G$ (identified with the product $G^{\abs{\mathcal W}}$), as the product
measure, $\nu_{\alpha,{\mathcal W}}=\nu_\alpha^{\otimes {\abs{\mathcal W}}}$.
The net $\left( \nu_{\alpha,{\mathcal W}}\right)$ has the desired properties:
concentration follows e.g. from \cite{Ta}, p. 76 and Prop. 2.1.1, while
the condition ($\ast$) is checked directly.
\end{proof}

In the case where $G$ is a compact group, the theorem was established by
Glasner \cite{Gl} and (independently, unpublished) Furstenberg and Weiss.
The above version is Thm. 2.2 in \cite{P5}. It would be interesting to know
if the result remains true for an arbitrary (non locally compact) amenable
topological group $G$. (Here of course one needs to replace measurable
maps with either Borel measurable ones, or else strongly measurable, that is,
measurable in the sense of Bourbaki, cf. \cite{Gaal}, p. 357.)

\begin{remark} Notice that, for $G=\R$, the topological group
$L^0(X)=L^0(X,\R)$ is extremely amenable, but not a L\'evy group, nor indeed
a generalized L\'evy group
in the sense of our Definition \ref{glg}, because it contains no non-trivial
compact subgroups.
%
\label{eanonl}
\end{remark}

Let us finally state a simple and useful fact that belongs to the folklore.

\begin{lemma}
\label{dense}
Let $G$ be a topological group, and let $\mathcal G$ be a family of
topological subgroups of $G$ with the property: if $H,K\in {\mathcal G}$,
then $H\cup K\subseteq L$ for some $L\in {\mathcal G}$. 
Suppose that $\cup\{G\colon G\in{\mathcal G}\}$ is everywhere dense in $G$ and that every
$H\in {\mathcal G}$ is contained in a suitable
extremely amenable subgroup of $G$. Then $G$ is extremely amenable.
\end{lemma} 

\begin{proof}
Let $G$ act continuously on a compact space $X$. For every $H\in {\mathcal G}$
denote by $F_H$ the set of all $H$-fixed points in $X$. As $H$ is contained
in a suitable extremely amenable subgroup of $G$, 
$F_H$ is non-empty and compact. As $\mathcal G$ is
upward directed, the family $\{F_H\colon H\in {\mathcal G}\}$ is centered, 
that is, for every finite subset
$H_1,\ldots,H_k\in {\mathcal G}$, one has $\cap_{i=1}^n F_{H_i}\neq\varnothing$.
Hence, $\cap_{H\in{\mathcal G}} F_H\neq\varnothing$,
and every element of this intersection is a fixed point for $G$.
\end{proof}

\begin{corollary}
Let $G$ be a topological group containing a family of locally
compact amenable subgroups directed
by inclusion and having an everywhere dense union. 
Let $X$ be a non-atomic Lebesgue measure space. Then 
the group $L^0(X,G)$ of all (equivalence classes) of Borel
measurable maps from $X$ to $G$, 
equipped with the topology of convergence in measure, is extremely amenable. 
\label{big}
\end{corollary}

\begin{proof} Let $\mathcal K$ be a family of locally compact
amenable subgroups of $G$ as in the
assumptions of the theorem. The group $L^0(X,\cup{\mathcal{K}})$ is everywhere
dense in $L^0(X,G)$, and it is extremely amenable by Lemma \ref{dense}, because
each of the topological subgroups of the form $L^0(X,K)$, 
$K\in{\mathcal{K}}$ is extremely amenable by Theorem \ref{gener}.
\end{proof} 

\section{Unitary groups of operator algebras}

\subsection{A characterization of approximately finite dimensional von
Neumann algebras}
Let $\H$ be a Hilbert space, let $M\subset \B(\H)$ be a von Neumann algebra,
and let $M_\ast$ be the predual of $M$. Recall that the 
{\it strong topology} on $M$, or the $s(M,M_\ast)$-topology, is the
topology on $M$ induced by the family of semi-norms
\[x\in M\mapsto \norm x_{\varphi} =\varphi(x^\ast x)^{1/2},~~ 
\varphi\in M_\ast^+,\]
and the {\em weak topology} on $M$ is the $\sigma(M,M_\ast)$-topology.
Then (see for example \cite{Tak}, Rem. II.4.10), the weak and the strong 
topologies coincide on the unitary group of $M$.

Let us mention two technical results that we will need.

\begin{lemma} 
\label{product}
Let $A$ be an index set and let $(M_\alpha)_{\alpha\in A}$ be
a family of \vNa s. 
If $M$ denotes the direct sum $\oplus_{\alpha\in A} M_\alpha$,
then $U(M)$, endowed with the $s(M,M_\ast)$-topology, is isomorphic
(as a topological group) to $\prod_{\alpha\in A} U(M_\alpha)$, where each
$U(M_\alpha)$ is endowed with the $s(M_\alpha,(M_\alpha)_\ast)$-topology
and $\prod_{\alpha\in A} U(M_\alpha)$ with the product topology.
\end{lemma}

\begin{proof}
As
\[M_\ast = \{(\varphi_\alpha)_{\alpha\in A} \colon \varphi_\alpha \in
(M_\alpha)_\ast,~~ \sum_{\alpha\in A}\norm{\varphi_\alpha}<\infty\},\]
the proof is clear.
\end{proof}

Let $N\subset\B(\H)$ be a von Neumann factor acting on a separable Hilbert
space $\H$ and let $M=L^\infty (X,\mu)\otimes N$, where $(X,\mu)$ is a 
Lebesgue space.

\begin{lemma} 
\label{l:u(m)}
With the above notation, the unitary group $U(M)$, with the
$s(M,M_\ast)$-topology, is isomorphic to the group $L^0(X,\mu; U(N))$ of all
measurable maps from $X$ to $U(N)$, equipped with the topology of convergence
in measure.
\end{lemma}

\begin{proof} Recall (\cite{Tak}, Thm. IV.7.17) that 
$(L^\infty(X,\mu)\otimes N)_\ast=L^1_{N_\ast}(X,\mu)$ and that 
${\mathrm{span}}\{f\otimes \varphi\colon f\in L^1(X,\mu),~\varphi\in N_\ast\}$
is norm-dense in $L^1_{N_x}(X,\mu)$. 

Let $(u_n)_{n\geq 1}$ be a sequence in $L^0(X,\mu; U(N))$ converging in measure
to $u\in L^0(X,\mu; U(N))$. Then $u_n\to u$ weakly, as for any $\varphi\in N_\ast$
and $f\in L^1(X,\mu)$, $f\otimes \varphi(u_n)\to f\otimes \varphi(u)$.
Indeed, if for $x\in X$, $g_n(x)=f(x)\varphi(u_n(x))$ and
$g(x) = f(x) \varphi(u(x))$, then $(g_n)_{n\geq 1}$, $g\in L^1(X,\mu)$
and
\[\abs{g_n(x)}\leq \abs{f(x)}\cdot \norm\varphi\mbox{ and }
\abs{g(x)}\leq \abs{f(x)}\cdot \norm\varphi,~~\mbox{$\mu$-a.e.}\]
Hence 
\begin{eqnarray*}
\abs{f\otimes\varphi(u_n) - f\otimes\varphi(u)} &\leq& 
\int\abs{f(x)}\abs{\varphi(u_n(x))-\varphi(u(x))}~d\mu(x) \\
&\leq & \int \abs{g_n(x)-g(x)}~d\mu(x) \underset{n\to\infty}{\rightarrow} 0.
\end{eqnarray*}
Conversely, let $\varphi\in N_\ast$ be a state on $N$ and let $\e>0$.
Let $(u_n)_{n\geq 1}$ and $u$ belong to $L^1(X,\mu; U(N))$. Then for $x\in X$,
\[\abs{\varphi(u_n(x))-\varphi(u(x))}=\abs{\varphi(u_n(x)-u(x))}
\leq \norm{u_n(x)-u(x)}_{\varphi}.\]
Hence (see \cite{cohn}, Prop. 3.1.5),
\[\mu\left(\{x\in X\colon \abs{\varphi(u_n(x))-\varphi(u(x))}\geq \e \} \right)
\leq
\mu\left(\{x\in X\colon \norm{u_n(x)-u(x)}^2_{\varphi}\geq \e^2\} \right).\]
If $u_n\to u$ strongly, then
\[\norm{u_n-u}_{f\otimes\varphi} =\left(\int \norm{u_n(x)-u(x)}^2_{\varphi}
~d\mu(x)\right)^{1/2}\underset{n\to\infty}{\rightarrow} 0.\]
By Tchebycheff's inequality, $(u_n)_{n\geq 1}$ converges in measure to $u$.
\end{proof}

The main result of this section is the following.

\begin{theorem}
If $M$ is an injective von Neumann algebra, then its unitary group
$U(M)$ with the $s(M,M_{\ast})$-topology is the direct product of
a compact group and an extremely amenable group.
\label{direct}
\end{theorem}

Let $M$ be an injective \vNa. By \cite{elliott}, Thm. 4, $M$ 
is generated by an 
upwards directed collection $(M_\alpha)_{\alpha\in A}$ of countably
generated injective von Neumann subalgebras. 
By \cite{dlH1}, Lemma 4, $U(M)$ is the
$s(M,M_\ast)$-closure of the upward directed collection of groups $U(M_\alpha)$.
(The separability condition in \cite{dlH1} is not essential, as 
remarked by Haagerup in \cite{haagerup}.) 
A countably generated injective \vNa\ is a
direct sum of injective \vNa s with separable predual. By Lemma \ref{product},
we can therefore assume in the proof of Theorem \ref{direct} that
$M_\ast$ is separable. 

Moreover, $M$ can be decomposed as a direct sum of finite \vNa s of type I
and type II, and of a properly infinite \vNa. Therefore 
Theorem \ref{direct} will follow directly from Proposition \ref{3.3},
Corollary \ref{3.5} and Proposition \ref{3.6}.

%
%
%
%

\begin{proposition} If $M$ is a finite \vNa\ of type $I$, then $U(M)$, with
the $s(M,M_\ast)$-topology, is the direct product of a family of compact groups and
extremely amenable groups.
\label{3.3}
\end{proposition}

\begin{proof} By the structure theorem for finite type $I$ von Neumann algebras
(see \cite{Tak}, Theorem V.1.27), $M$ is isomorphic to a direct sum of
algebras of the form $L^\infty(\Gamma,\mu)\otimes \B(\H)$, where 
$\dim\H<\infty$ and $\Gamma$ is a locally compact space, with a positive
(finite) Radon measure $\mu$. 

If $\mu$ is atomic, then
\[U(L^\infty(\Gamma,\mu)\otimes \B(\H))\cong \prod_{\supp\mu} U(\H)\]
and therefore is a compact group.

If $\mu$ is non-atomic, then the unitary 
group $U(L^\infty(\Gamma,\mu)\otimes \B(\H))$ is canonically isomorphic to 
$L^\infty (\Gamma,\mu; U(\H))$. 
It follows that $U(L^\infty(\Gamma,\mu)\otimes \B(\H))$ is extremely amenable by
Lemma \ref{l:u(m)}.
By Lemma \ref{product} and Corollary \ref{big}, the Proposition is proved.
\end{proof}

\begin{proposition}
\label{increasing}
Let $M$ be a \vNa\ and let $(N_k)_{k\geq 1}$ be an increasing sequence of
finite-dimensional subfactors of $M$ such that 
$N_\infty=\cup_{k\geq 1}N_k$ is $s(M,M_\ast)$-dense in $M$. 
Then $U(M)$, with the $s(M,M_\ast)$-topology, is a L\'evy group.
\end{proposition}

\begin{proof} For $k\geq 1$, let $(e_{i,j}^k)_{1\leq i,j\leq n_k}$ be a
system of matrix units of $N_k$. If $u\in U(N_k)$, let $(u_{ij})\in U(n_k,\C)$
denote the matrix defined by $u=\sum_{i,j=1}^k u_{ij} e_{ij}^k$. Let
$SU(N_k)=\{u\in U(N_k)\colon \det(u_{ij})=1\}$. 

Let $\varphi\in M_\ast$ be a faithful state and let $\e>0$ be given.
As $\varphi(1)=\sum_{i=1}^{n_k}\varphi(e_{ii}^k)$, there exists $K$ such
that for $k\geq K$, $\varphi(e_{ii}^k)\leq \e/4$ for some $1\leq i\leq n_k$.

Let $u\in U(N)$. By \cite{dlH1}, Lemma 4, there exist $k\geq K$ 
and $u_k\in U(N_k)$ such that $\norm{u-u_k}_{\varphi}\leq \e/2$. 
Let $1\leq i\leq n_k$ be such that $\varphi(e_{ii}^k)\leq \e/4$ and
$v_k = 1-e_{ii}^k+\overline{\det(u_{ij}^k)}e^k_{ii}\in N_k$. Then
$u_kv_k\in SU(N_k)$ and 
\begin{eqnarray*}
\norm{u-u_kv_k}_{\varphi} &\leq & 
\norm{u-u_k}_{\varphi} + \norm{u_k(1-v_k)}_{\varphi} \\
&\leq & \frac \e 2 + \norm{1-v_k}_{\varphi} \\
&\leq& \frac \e 2 + \lvert 1 - \det(u_{ij}^k)\rvert \varphi(e_{ii}^k) \\
&\leq& \frac \e 2 +2\frac \e 4 = \e.
\end{eqnarray*}
Hence, $\cup_{k\geq 1} SU(N_k)$ is $s(M,M_\ast)$-dense in $U(M)$.

As $N_k\cong M_{n_k}(\C)$, there exists $h_{\varphi} \in M_{n_k}(\C)$,
$h_{\varphi}> 0$ such that for 
$x=\sum_{i, j = 1}^{n_k} x_{i,j} e^k_{i,j} \in N_k$, 
$\varphi(x) =\Tr (h_{\varphi}(x_{i,j}))$,
where $\Tr$ denotes the trace on $M_{n_k}(\C)$. As 
$\varphi(1) =\Tr(h_{\varphi})=1$,
$\norm{h_{\varphi}}\leq 1$ and for $x\in N_k$,
\[\norm x_{\varphi}\leq \norm{(x_{ij})}_2,\]
where $\norm\cdot_2$ denotes the Hilbert--Schmidt norm. 
By \cite{GrM}, Example 3.4 $(SU(N_k))_{k\geq 1}$ forms a L\'evy family with regard
to the Hilbert-Schmidt metric, and consequently the 
$\norm\cdot_\varphi$-metric, which ends the proof of Proposition.
\end{proof} 

Let $M$ be a properly infinite injective \vNa\ with separable predual.
By \cite{elliott}, there exists an increasing sequence 
$(N_k)_{k\geq 1}$ of finite dimensional subfactors whose union is
$s(M,M_\ast)$-dense in $M$. By Proposition \ref{increasing}, we get

\begin{corollary}
\label{3.5}
Let $M$ be a properly infinite injective \vNa\ with separable predual.
Then the unitary group 
$U(M)$, with the $s(M,M_\ast)$-topology, is extremely amenable. \qed
\end{corollary}

Let $M$ be a finite continuous injective \vNa\ with separable predual. By
considering the standard representation of $M$, we can assume that $M$ acts
on a separable Hilbert space. Then by reduction theory and uniqueness of the 
injective factor of type $II_1$, $M$ is isomorphic to
$L^\infty(X,\mu)\otimes R$, where $X$ is a standard Borel space, with a 
finite measure $\mu$ and $R$ is the hyperfinite factor of type $II_1$.

By Proposition \ref{increasing} and Corollary \ref{everylevy}, $U(R)$, 
endowed with the $s(R,R_\ast)$-topology,
is extremely amenable. By the same argument as in the proof of \ref{3.3},
if $\mu$ is an atomic measure, then $U(M)$
is isomorphic to $\prod_{\supp\mu} U(R)$, equipped
with the product $s(R,R_\ast)$-topology. If $\mu$ is a non-atomic measure, then
$U(M)$, with the $s(M,M_\ast)$-topology, is isomorphic to the group
of all measurable maps from $X$ to $U(R)$, equipped with the topology of
convergence in measure. Then by Lemma \ref{product} and Corollary \ref{big}
we get 

\begin{proposition}
\label{3.6}
If $M$ is a finite continuous injective
\vNa\ with separable predual, then $U(M)$,
endowed with the $s(M,M_\ast)$-topology, is extremely amenable. \qed
\end{proposition}

\subsection{A characterization of nuclear $C^\ast$-algebras}

Recall that a continuous action of a topological group $G$ on a compact
space $X$ is {\it proximal} if for every two points $x,y$ and every entourage
of the diagonal $V\subseteq X\times X$ there is a $g\in G$ with $(gx,gy)\in V$.
A topological group $G$ is called {\it strongly amenable} if
every continuous proximal action of $G$ on a compact space has a fixed point
\cite{Gl1}. For instance, every compact group is strongly amenable simply 
because it has no non-trivial proximal actions on compacta. Every extremely
amenable group is strongly amenable for obvious reasons. The class of
strongly amenable groups is closed under direct products. 

By the {\it completion} $\hat G$ 
of a topological group $G$ we mean, as usual, the
completion of $G$ with regard to the two-sided uniform structure, that is,
the supremum of the left and the right uniform structures on $G$. This is
again a topological group.
(See, e.g., Chapter 10 in \cite{RD}.)

It is a standard fact, easily proved, that the completion, $\hat G$, 
of a topological group $G$ is amenable (strongly amenable,
extremely amenable) if and only if $G$ is amenable (resp. strongly amenable,
extremely amenable). 

The following corollary of our results 
strengthens a result due to Paterson \cite{Pat} (Thm. 2).

\begin{theorem} A $C^\ast$-algebra $A$ is nuclear if and only if  
its unitary group $U(A)$, equipped with the
topology $s(A,A^\ast)$, is strongly amenable. 
\end{theorem} 

\begin{proof} Denote by $M$ the universal von Neumann envelope of $A$.
As is well known, $M$ can be obtained as the completion 
of the $C^\ast$-algebra $A$
with regard to the additive uniform structure determined by the
relative weak topology. (\cite{dixmier}, Corollaire 12.1.3.) 
As a consequence, one can prove that the 
completion, $\widehat{U(A)}$, of the unitary group 
$U(A)$, equipped with the relative weak
topology, is isomorphic to the unitary group of the von Neumann envelope
$M$ of $A$, equipped with the $s(M,M_\ast)$-topology. 
\smallskip

Necessity $(\Rightarrow)$: if $A$ is nuclear, then the von Neumann envelope
$M$ is approximately finite dimensional, and therefore, by Theorem
\ref{direct}, the unitary group
$U(M)$ with the topology $s(M,A^\ast)$ is the product of a compact
group and an extremely amenable group. In particular, $U(M)$ is strongly
amenable, and so is its everywhere dense topological subgroup $U(A)$.
\smallskip

Sufficiency $(\Leftarrow)$: follows from de la Harpe's result in
\cite{dlH1}, as noted by Paterson \cite{Pat}.
If the topological group $U(A)$ is strongly
amenable, then it is in particular amenable, and so is the topological
group completion, $U(M)$. It follows by de la Harpe's result
 that $M$ is approximately finite-dimensional.
\end{proof}

\subsection{Automorphism groups of von Neumann algebras\label{agvna}}

If $M$ is a \vNa, let $\Aut(M)$ denote, as usual, the group of all
$\ast$-automorphisms of $M$. We consider on $\Aut(M)$ the topology of norm
pointwise convergence for the action of $\Aut(M)$ on $M_\ast$, given by
\[\alpha\in\Aut(M),~~\varphi\in M_\ast \mapsto \varphi\circ \alpha^{-1}
\in M_\ast.\]
With this topology, called the u-{\it topology}, $\Aut(M)$ is a topological
group and it is a Polish group if $M_\ast$ is separable.

If $u\in U(M)$, let $\Ad u$ denote the inner automorphism of $M$ given by
$\Ad u(x) =uxu^\ast$, $x\in M$, and let $\Inn(M)$ be the subgroup of all
inner automorphisms of $M$. 

If $U(M)$ is endowed with the $s(M,M_\ast)$-topology and $\Inn(M)$ with the
$u$-topology, the canonical surjection $U(M)\to\Inn(M)$ is continuous, as for
$\varphi\in M_\ast$, $u\in U(M)$
\[\norm{\varphi\circ\Ad u -\varphi}\leq 2\norm\varphi\norm{u-1}_{\varphi}.\]
By Theorem \ref{direct}, we therefore get:

\begin{proposition}
If $M$ is an injective continuous \vNa, then $\Inn(M)$, with the topology of
pointwise convergence in norm on $M_\ast$, is extremely amenable. \qed
\label{pcn}
\end{proposition}

\section{Groups of measure-preserving transformations}

\subsection{\label{twotopologies}Weak and uniform topologies}

For a standard non-atomic Lebesgue measure space $X=(X,\mu)$ denote by
$\Aut(X,\mu)$ the group of (equivalence classes of)
invertible measure-preserving transformations of
$X$, and by $\Aut^\ast (X,\mu)$
the group of all invertible measurable and
non-singular transformations of $X$. 

The {\it weak topology} on $\Aut^\ast (X,\mu)$ (also called the
{\em coarse topology}) is 
induced by the strong operator topology on the isometry group of
$L^p(X)$ under the quasi-regular representation of 
the group $\Aut^\ast (X,\mu)$ in the Banach space $L^p(X)$, where 
$1\leq p<\infty$ is any. (Cf. \cite{CK}, Theorem 8.) 

The weak topology makes $\Aut^\ast (X,\mu)$ into a Polish topological group,
while $\Aut\,( X,\mu)$ is a closed (therefore also Polish) subgroup. 
(See \cite{K}, 17.46.) Notice that the weak topology makes perfect sense
not just for a finite measure $\mu$, but 
for sigma-finite one as well.

%

Define a left-invariant metric on $\Aut^\ast (X,\mu)$ as follows:
\[d_{unif}(\tau,\sigma)=\mu\{x\in X\colon \tau(x)\neq\sigma(x)\}.\]
The topology induced on $\Aut^\ast(X,\mu)$ by the metric $d$ is a group topology,
%
known as the {\it uniform topology.} It is
strictly finer than the weak
topology (\cite{IT}, Prop. 3). For example, the uniform topology 
makes $\Aut^\ast (X,\mu)$ into a (path-connected) non-separable
group. When dealing with the uniform topology, we will always assume $\mu$
to be a finite measure.

Both the weak topology and the uniform topology 
only depend on the equivalence class of the measure $\mu$, rather
than on $\mu$ itself. (See \cite{CP}, Rem. 3, p. 373.)

Sometimes we will indicate the uniform topology by the subscript $u$,
as in $\Aut(X,\mu)_u$. Similarly, for the weak topology the subscript $w$
will be used.

Note that the restriction of $d$ to $\Aut(X,\mu)$ is bi-invariant, and so the
topological group $\Aut(X,\mu)_u$ is SIN (= has small invariant neighbourhoods),
that is, the left and the right uniform structures coincide.


\subsection{Extreme amenability of $\Aut(X,\mu)$ with the weak topology}
The following theorem (cf. e.g. \cite{Hal},
pp. 65--68), belonging to the Kakutani--Rokhlin circle of results,
is well-known. 
Identify the symmetric group ${\mathfrak{S}}_n$ 
with the subgroup of all measure-preserving 
automorphisms of ${\mathbb I}=[0,1]$ with the Lebesgue measure $\lambda$
mapping each dyadic interval of rank $n$ onto a dyadic interval of rank $n$
via a translation (the interval exchange transformations).

\begin{theorem}[Weak Approximation Theorem]
The union of the subgroups ${\mathfrak{S}}_n$, $n\in\N$
is everywhere dense in $\Aut({\mathbb I},\lambda)$ with respect to the weak (=strong)
topology. \qed
\label{wat}
\end{theorem}

Under the above embedding ${\mathfrak{S}}_n\hookrightarrow \Aut(\I)$, 
the normalized Hamming distance $d_n$ on ${\mathfrak{S}}_n$
(Example \ref{maurey}) coincides with the restriction of the metric $d$. 

\begin{theorem} 
The group $\Aut(X,\mu)_w$ of all measure-preserving 
automorphisms of a standard non-atomic finite or sigma-finite 
measure space, equipped with 
the weak topology, is a L\'evy group and therefore extremely amenable.
\label{isextre}
\end{theorem}

\begin{proof} For $\mu$ finite, the group $\Aut(X,\mu)$ is topologically
isomorphic with the group $\Aut(\I,\lambda)$, where 
$\I=[0,1]$ is the unit interval with the Lebesgue measure.
According to Maurey's theorem \cite{Ma}, 
the permutation groups ${\mathfrak{S}}_n$ of finite rank $n$, equipped
with the normalized counting measure and the normalized Hamming distance, form
a L\'evy family (even with regard to the uniform topology). By 
Theorem \ref{wat}, $\Aut(X,\mu)_w$
forms a L\'evy group, and Corollary \ref{everylevy} applies.

For $\mu$ sigma-finite it is enough to prove the result for the real line
$\R$, equipped with the Lebesgue measure $\lambda$. 
For every $n\in\N$, let $G_n$ denote the group of interval exchange 
transformations of the system of 
$8^{n+1}$ intervals of length $16^{-n-1}$ each,
covering the interval $[-2^{n},2^{n}]$ of length $2^{n+1}$. 
A natural group isomorphism from
the permutation group ${\mathfrak{S}}_{8^{n+1}}$ of rank $8^{n+1}$ 
equipped with the normalized
Hamming distance onto $G_n$, equipped with the uniform metric, 
is Lipschitz with $L=2^{n+1}$. By Lemma \ref{above}
and Maurey's result (Example \ref{maurey}), the family $(G_n)$ is a L\'evy
family, with 
\begin{eqnarray*}
\alpha_{G_n}(\e) &\leq & \alpha_{{\mathfrak{S}}_{8^{n+1}}}(\e/2^{n+1}) \\
&\leq & \exp\left(-\frac{\e^2}{4^{n+1}}\cdot\frac{8^{n+1}}{32} \right) \\
&=& \exp\left(-\frac{\e^22^{n+1}}{32} \right)\to 0\mbox{ as }n\to\infty.
\end{eqnarray*}
Since the groups $(G_n)$ form an increasing chain with everywhere dense
union in $\Aut(\R,\lambda)$, we are done.
\end{proof}

\subsection{Non-amenability of $\Aut(X,\mu)$ with uniform topology}
Let $\pi$ be a unitary representation of a group $G$ in a Hilbert space $\H$,
that is, $\H$ is a unitary $G$-module. One says that $\pi$ (or
the $G$-module $\H$) is {\it amenable in the sense of Bekka} \cite{Bek1}
if there exists
a state, $\phi$, on the von Neumann algebra ${\mathcal B}(\H)$ that is
invariant under the action of $G$ by conjugations: for all $T\in {\mathcal B}(\H)$
and all $g\in G$,
\[\phi(T) = \phi(\pi^\ast_g T\pi_g).\]
Denote by $\s$ the unit sphere in the space $\H$.  Notice that $G$ acts on
the space $\C^{\s}$ of all functions
$\s\to\C$ by 
\[^g\! f(\xi):=f(\pi^\ast_g(\xi))\equiv f(\pi_{g^{-1}}(\xi)),\]
where $g\in G$, $f\in \C^{\s}$, and $\xi\in\s$.
This action leaves invariant the space $\UCB(\s)$ of all
uniformly continuous bounded complex-valued functions on the sphere. 
The following was proved by one of the present authors \cite{P4}.

\begin{theorem} 
A unitary representation $\pi$ of a group $G$ is amenable if and only if
there exists a $G$-invariant mean on the space $\UCB(\s)$. \qed
\end{theorem}

Let $\xi\in \s$. 
For every (complex-valued) 
function $f$ on $\s$ denote by 
$\tilde f$ a function on $G$ defined as follows:
\[\tilde f(g):= f(\pi_g\xi),\]
where $g\in G$. 
If $f\in \ell^\infty(\s)$, then $\tilde f\in \ell^\infty(G)$. 
The latter is a Banach $G$-module with respect to
the left regular representation of $G$: 
\[^g\!\varphi(x) := \varphi(g^{-1}x).\]
The mapping $f\mapsto \tilde f$ from $\ell^\infty(\s)$   to $\ell^\infty(G)$ is
$G$-equivariant, that is, commutes with the action of $G$:
\[\widetilde{^g\!f}(h) = ^g\!f(\pi_h(\xi)) 
= f(\pi^\ast_g\pi_h(\xi)) 
= f(\pi_{g^{-1}h}(\xi)) 
= \tilde f(g^{-1}h) 
= ^g\!(\tilde f)(h).\]
%

Now assume that $G$ is a topological group and the representation $\pi$
is strongly continuous. 

\begin{lemma} If $f\in\UCB(\s)$, then $\tilde f$ is left uniformly 
continuous on $G$.
\end{lemma}

\begin{proof}
Let $\e>0$. Find a $\delta>0$ such that, whenever $\norm{\xi-\zeta}<\delta$, one
has $\abs{f(\xi)-f(\zeta)}<\e$. 
Due to the strong continuity of $\pi$, there is a neighbourhood of
the identity, $V$, in $G$ such that, whenever $g\in V$, one has
$\norm{\xi-\pi_g\xi}<\delta$. If now $g,h\in G$ are such that $g^{-1}h\in V$, then
$\abs{\pi_{g^{-1}h}(\xi)-\xi}<\delta$, that is, 
$\abs{\pi_{h}(\xi)-\pi_g(\xi)}<\delta$, and
$\abs{\tilde f(g) - \tilde f(h)} = \abs{f(\pi_g(\xi)) - f (\pi_h(\xi))}
< \e.$
\end{proof} 

Denote by $\LUCB(G)$ (respectively, $\RUCB(G)$)
the $C^\ast$-algebra of all bounded left (respectively, right) uniformly 
continuous complex-valued functions on $G$.
It is now easy to verify that the mapping $f\mapsto \tilde f$ is a positive
linear operator from $\UCB(\s)$ to $\LUCB(G)$, sending $1$ to $1$ and commuting
with the action of $G$. Therefore, if there is an invariant mean on $\LUCB(G)$,
by composing it with the mapping $f\mapsto \tilde f$, one gets 
an invariant mean on $\UCB(\s)$. We obtain:

\begin{proposition}
If a topological group $G$ admits a left-invariant
mean on the space $\LUCB(G)$ of all left uniformly continuous bounded functions
on $G$, then every strongly continuous unitary representation of $G$ is
amenable. 
\label{pul} \qed
\end{proposition} 

Recall that a topological group is called {\em amenable} if there exists a
left-invariant mean on the space $\RUCB(G)$. 
If $G$ is a SIN group, that is, the left and the right 
uniformities on $G$ coincide, then $\LUCB(G)=\RUCB(G)$ and
we obtain the following.

\begin{corollary} 
\label{c:pul}
Every strongly continuous unitary representation of an
amenable SIN group is amenable. \qed
\end{corollary} 

\begin{remark} In general, the condition that $G$ is SIN cannot be dropped here,
unless we presume that $G$ is locally compact. For
instance, the infinite unitary group $U(\ell_2)$ with the strong operator
topology is (extremely) amenable, while the standard unitary representation
of $U(\ell_2)$ is non-amenable (because this would clearly imply, for instance,
amenability of every unitary representation of every countable group).
\end{remark}

\begin{theorem}
The topological group $\Aut(X,\mu)_u$ is non-amenable.
\end{theorem}

\begin{proof} Take as $V=L^2_0(X)$ the closed unitary
$G$-submodule of $L^2(X)$ consisting of all functions $f\in L^2(X)$
with
\[\int_X f(x)\, d\mu(x)=0.\]
According to Corollary \ref{c:pul}, if $\Aut(X,\mu)_u$ were amenable,
then the regular representation $\gamma$ of $\Aut(X,\mu)_u$ in 
$L^2_0(X)$ would be amenable. (Recall that with regard to the uniform topology
the group $\Aut(X,\mu)$ is a SIN group.)
In particular, the
restriction of $\gamma$ to any subgroup $G$ of $\Aut(X,\mu)_u$ 
would be amenable as well.
It is, therefore, enough to discover a group $G$ acting on a
standard non-atomic Lebesgue space $X$ by measure-preserving
transformations in such a manner that
the associated regular representation of $G$ in $L^2_0(X)$ is
not amenable in the sense of Bekka.

The following is inspired by Bekka's work \cite{Bek2}.
Let $G$ be a semisimple real Lie group with finite centre and
without compact factors, having property $(T)$. 
(Cf. \cite{Bek2}; e.g. $SL_3(\R)$.) 
Let $\Gamma$ denote a lattice in $G$, that is, a discrete subgroup
such that the Haar measure induces an invariant finite measure on
$G/\Gamma$. Set $X=G/\Gamma$. 

Assume that the standard representation $\pi$ of $G$ in
$L^2_0(X)=L^2_0(G/\Gamma)$ is amenable. Since $G$ is a Kazhdan group,
according to a result by Bekka \cite{Bek1}, 
$\pi$ must have a finite-dimensional
subrepresentation. (Later it was shown \cite{BV} that this property
characterizes groups with property $(T)$.)

Note that $L^2_0(X)$ has nontrivial $G$-invariant
vectors: since the action of $G$ on $X$ is transitive, such vectors
must be constant functions. Therefore, a finite-dimensional 
subrepresentation of $\pi$ must be non-trivial. 
But this would result in a non-trivial continuous
group homomorphism from $G$ to a finite-dimensional unitary group,
which is impossible by assumption.
\end{proof}

\section{Full groups of amenable equivalence relations}

\subsection{The full group\label{fullgr}}
Let $\mathcal R$ be a discrete measured equivalence relation, i.e.
$\mathcal R$ is an equivalence relation with countable equivalence classes
on a standard Borel space $S$, the graph $\mathcal R\subset S\times S$ is
Borel, and $\mu$ is a quasi-invariant probability measure on $S$. 

The full group $[{\mathcal R}]$ is defined as the group of all bimeasurable
transformations $\sigma$ of $(S,\mu)$ with
\[(s,\sigma (s))\in{\mathcal{R}},~~\mbox{$\mu$-a.e.}\]

We will consider on $[{\mathcal R}]$ the
uniform topology, determined by the measure $\mu$. With this topology,
$[{\mathcal R}]$ is a Polish group \cite{HO},\cite{danilenko}. 

Following  Zimmer \cite{zimmer-inv},  let $G$ be a countable discrete group and 
let $(S,\mu)$ be a right free $G$-space.
Then an element $\sigma$ of the {\it full group} $[G]$ of $G$ is
determined by a measurable partition with the property
\[S = \coprod_{g\in G} S_g = \coprod_{g\in G}S_g g,\]
in such a way that $s\sigma = sg$ if $s\in S_g$.
In particular, $G$ is a subgroup of $[G]$, every element $h\in G$ being determined
by the partition
\[S_g = \begin{cases} S, &\mbox{if $g=h$,} \\
\emptyset, &\mbox{otherwise.}
\end{cases}\]


\begin{remark} To the measured dynamical system $(S,\mu,G)$ is associated a
\vNa, the so-called Murray--von Neumann group
measure space construction. One has the following unitary representations
$U,V$ of $G$ and isometric $\ast$-representations $M$ and $N$ of $L^\infty(S)$ 
on $L^2(S\times G)$, given for $g\in G$, $f\in L^\infty(S,\mu)$, by
\[U_g\xi(x,h) = r(x,g)^{1/2}\xi(xg,hg),~~V_g\xi(x,h) = \xi(x,g^{-1}h),\]
\[M_f\xi(x,h) = f(x)\xi (x,h),~~ N_f\xi(x,h) = f(xh^{-1})\xi(x,h),\]
where $r(x,g)=d\mu(xg)/d\mu(x)$ is the 
Radon--Nikodym cocycle of the action of $G$ on $(S,\mu)$.

The \vNa\ generated by the operators $V_g$ and $N_f$ is denoted by $L$,
while that generated by the operators $U_g$ and $M_f$ is denoted by $R$.
\end{remark}

If $J\in U(L^2(S\times G))$ is given by
\[J\xi (x,g) = f(xg^{-1},g^{-1})~r(x,g^{-1})^{1/2},~~(x,g)\in S\times G,\]
then $J^2= 1$, and $JLJ = R$.

To any $\sigma\in [G]$, defined by the partitions $(S_g)_{g\in G}$ and
$(S_gg)_{g\in G}$ of $S$, there are naturally associated unitaries
$U_\sigma\in R$ and $V_\sigma\in L$ given by
\[U_\sigma =\sum_{g\in G}M_{f_g}U_g\in R\mbox{ and }
V_\sigma = \sum_{g\in G}N_{f_g}V_g\in L,\]
where $f_g=\chi_{S_g}$.

If $U([G])\subset U(R)$ (resp. $V([G])\subset U(L)$) is endowed with the
strong operator topology and $[G]$ with the uniform topology, it is easy to
check that $U$ (resp. $V$) is a topological isomorphism from $[G]$ onto its image. 

\subsection{\label{eafg}} In this subsection, we will show that if  
$\mathcal R$ is an amenable
equivalence relation, then the full group $[{\mathcal R}]$ 
with the uniform topology is extremely amenable. For a definition of an
amenable equivalence relation, see \cite{zimmer-inv}, Sect. 3.

By \cite{CFW}, we can assume that $S$ is the compact abelian group
$\prod \{0,1\}$, the equivalence relation ${\mathcal R}$ is induced by
the action of the subgroup $G=\oplus \Z/2$, acting by addition: if
$g=(g_k)_{k\geq 1} \in G$, and $x= (x_k)_{k\geq 1}$, then 
\[xg = (x_k+g_k)_{k\geq 1}=x+g,\]
and $\mu$ is a $G$-ergodic non-atomic quasi-invariant measure on $S$. 

Let us introduce some notation we will use in this section.
For $n\geq 1$, let $S_n = \prod_{k=1}^n \{0,1\}$, and 
$S^n=\prod_{k\geq n+1}\{0,1\}$. 
For $y\in S$ we
will denote its tail by $y^n = (y_k)_{k\geq n+1}\in S^n$.
For $x\in S_n$ let
$C(x)$ denote the cylinder set $\{(y_k)_{k\geq 1}\in S\colon
y_k =x_k, 1\leq k \leq n\}$ of $S$. 
For every bijection $\sigma\in \Bij(S_n)$,
let $\sigma$ be the element of $[G]$ defined for $y=(x,y^n)\in C(x)$, $x\in S_n$,
by $\sigma(y) = (\sigma(x),y^n)$, and let $\Bij(n)$ denote the finite
subgroup of $[G]$ such transformations form. 

For $n\geq 1$, if $G_n = \{g = (g_k)_{k\geq 1}\in G\colon g_k =0~~
\forall k\geq n+1\}$, then $\Bij(n)\subset [G_n]$.

\begin{lemma} Keeping the above notation, $\cup_n \Bij(n)$ is uniformly
dense in $[G]$.
\label{approximation}
\end{lemma}

\begin{proof} Let $(S_g)_{g\in G}$ and $(S_gg)_{g\in G}$ be a measurable
partition of $S$ and let $\alpha\in [G]$ be given by
$\alpha(x) = xg$, $x\in S_g$. Let $0<\e<1$ be given and let $N_1$ be such that
$\mu\left(\coprod_{g\in G_{N_1}}S_g\right)\geq 1 - \e$. Let us fix 
$0<\delta<\frac{\e}{3\cdot 2^{2N_1}}$.

The probability measure $\nu =\frac 1{2^{N_1}}\sum_{g\in G_{N_1}}\mu\circ g$
is $G_{N_1}$-invariant and equivalent to $\mu$.

For $n\geq 1$, let $E_n$ denote the conditional expectation with respect to
$\nu$ and to the partition by the cylinders of $S_n$, i.e. for a measurable subset
$B\subset S$,
\[E_n(\chi_B) = \sum_{x\in S_n} r^n_x(B)\chi_{C(x)},\]
where 
\[r^n_x(B) = \frac{\nu(C(x)\cap B)}{\nu(C(x))}.\]
Let $N_2\geq N_1$ be such that for all $B=S_g$
or $B=S_gg$, for $g\in G_{N_1}$,
\begin{equation}
\label{5.1}
\norm{E_{N_2}(\chi_B)-\chi_B}_{1,\nu}=2 \sum_{x\in S_{N_2}} r^{N_2}_x(B)
\left(1-r^{N_2}_x(B)\right)\nu (C(x))<\delta^2 \end{equation}
If $A$ is a measurable subset of $S$, let $S_{N_2}^A$ denote
$\{x\in S_{N_2}\colon r^{N_2}_x(A) \geq 1-\delta\}$ and let
$\tilde A$ denote the disjoint union 
$\coprod_{x\in S^A_{N_2}}C(x)$. 
Then the sets 
$(\tilde S_g)_{g\in G_{N_1}}$ satisfy for $g,h\in G_{N_1}$ 
the following properties:
\smallskip
\begin{itemize}
\item{i)} $\norm{\chi_{S_g}-\chi_{\tilde S_g}}_{1,\nu}\leq 3\delta$,
\smallskip
\item{ii)} $\tilde S_g\cap \tilde S_h =\emptyset$, for 
$g\neq h$,
\smallskip
\item{iii)} $\tilde S_gg\cap \tilde S_hh=\emptyset$, for $g\neq h$.
\end{itemize}
Indeed, for $g\in G_{N_1}$ we have by definition of $\tilde S_g$:
\begin{eqnarray*}
\norm{E_{N_2}(\chi_{S_g})-\chi_{\tilde S_g}}_{1,\nu} &\leq & 
\delta \sum \left\{\nu(C(x)) \colon x\in S_{N_2},~ r^{N_2}_x(S_g)\leq\delta
\mbox{ or }r^{N_2}_x(S_g)\geq 1-\delta\right\} +\\[2mm]
&& \sum\left\{\nu(C(x)) \colon x\in S_{N_2},~ r^{N_2}_x(S_g)(1-r^{N_2}_x(S_g))
>\delta(1-\delta)\right\}.
\end{eqnarray*}
By (\ref{5.1}) we have
\begin{eqnarray*}
\norm{E_{N_2}(\chi_{S_g})-\chi_{\tilde S_g}}_{1,\nu}
&\leq& \delta + \frac 1{2\delta(1-\delta)} \norm{E_{N_2}(\chi_{S_g})-\chi_{S_g}}_{1,\nu} \\
&\leq& \delta + \frac{\delta}{2(1-\delta)} \leq \frac{5\delta}2.
\end{eqnarray*}
Then,
\begin{eqnarray*}
\norm{\chi_{S_g}-\chi_{\tilde S_g}}_{1,\nu} &\leq& 
\norm{E_{N_2}(\chi_B)-\chi_B}_{1,\nu}
+ \norm{E_{N_2}(\chi_{S_g})-\chi_{\tilde S_g}}_{1,\nu} \\ 
&\leq&
\delta^2+\frac{5\delta}2 \leq
3\delta,
\end{eqnarray*}
which proves (i).

For $x\in S_{N_2}$,
\[
1\geq \sum_{g\in G_{N_1}} \frac{\nu(C(x)\cap S_g)}{\nu(C(x))}
= \sum_{g\in G_{N_1}} r^{N_2}_x(S_g) \geq 0.
\]
Hence $S^{S_g}_{N_2}\cap S^{S_h}_{N_2} =\emptyset$ and
$\tilde S_g\cap \tilde S_h =\emptyset$, which proves (ii). 

For (iii), notice that for $x\in S_{N_2}$, $g\in G_{N_1}$,
\[r^{N_2}_{xg}(S_g) = \frac{\nu(S_g\cap C(xg))}{\nu(C(xg))} =
\frac{\nu((S_gg\cap C(x))g)}{\nu(C(xg))} = r^{N_2}_x(S_gg),
\]
as $\nu$ is $G_{N_1}$-invariant. Hence $S^{S_g}_{N_2}g = S^{S_gg}_{N_2}$
and
\[\tilde S_gg = \coprod_{x\in S^{S_g}_{N_2}}C(xg) = 
\coprod_{x\in S^{S_g}_{N_2}g}C(x)=\coprod_{x\in S^{S_gg}_{N_2}}C(x).\]
If $g,h\in G_{N_1}$, $g\neq h$, then $S_gg\cap S_hh=\emptyset$ and (iii)
holds.

We now define $\sigma\in \Bij(N_2)$ such that $d_\mu(\alpha,\sigma)<2\e$.
As $S^{S_g}_{N_2}g = S^{S_gg}_{N_2}$, for $g\in G_{N_1}$, 
$A=\left\{ x\in S_{N_2}\colon x\notin \coprod_{g\in G_{N_1}}
S^{S_g}_{N_2}\right\}$ and $B=
\left\{ x\in S_{N_2}\colon x\notin \coprod_{g\in G_{N_1}}
\left(S^{S_g}_{N_2}g \right)\right\}$ have the same cardinality. 
Let $\tilde\sigma$ be an arbitrary bijection from 
$A$ to $B$, and let
$\sigma$ denote the bijection of $S_{N_2}$ defined by
\[ \sigma(x) = \begin{cases}
xg, & \mbox{ if }x\in S^{S_g}_{N_2}, g\in G_{N_1}, \\
\tilde \sigma(x), &\mbox{ if }x\in A.
\end{cases}\]
We shall also denote by $\sigma$ the corresponding element of $\Bij(N_2)$.
Then 
\[\{x\in S\colon \sigma(x)\neq \alpha(x)\}\subseteq 
\left(S\setminus \coprod_{g\in G_{N_1}}S_g \right) \coprod
\left( \coprod_{g\in G_{N_1}}(S_g\setminus\tilde S_g) \right)
\]
and 
\[d_\mu(\alpha,\sigma) \leq \mu
\left(S\setminus \coprod_{g\in G_{N_1}}S_g \right) +
\sum_{g\in G_{N_1}}\mu (S_g\Delta \tilde S_g).
\]
As $\mu (S_g\Delta \tilde S_g)\leq 2^{N_1} \nu(S_g\Delta \tilde S_g) =
2^{N_1}\norm{\chi_{S_g}-\chi_{\tilde S_g}}_{1,\nu} \leq
2^{N_1	}3\delta$, then
$d_\mu(\alpha,\sigma)\leq \e+2^{2N_1}3\delta\leq 2\e$.
\end{proof}

As the size of the largest cylinder set determined by the first $n$ coordinates
clearly goes to zero as $n\to\infty$, Corollary \ref{largestatom} implies
that the family of groups $\Bij(n)$ equipped with the uniform metric and
the normalized counting measure is a L\'evy family. 
Applying Lemma \ref{approximation}, we deduce:

\begin{proposition} 
\label{hard}
Let $\mathcal R$ be a discrete measured equivalence relation, acting 
ergodically on a Lebesgue space $(S,\mu)$. If $\mathcal R$ is an
amenable equivalence relation, then the full group $[{\mathcal R}]$, 
endowed with the uniform
topology, is a L\'evy group (and in particular
extremely amenable). \qed
\end{proposition}

\begin{remark}
As pointed out to us by A.S. Kechris, extreme amenability of the full group
$[G]$ in Proposition \ref{hard} can be alternatively
deduced by applying an argument
similar to that employed by one of the present authors to give a new proof
of extreme amenability of the unitary group $U(\ell^2)_s$ in \cite{P5}, 
subsection 4.5.

Namely, for $n\geq 1$, let
$S^n = \prod_{k\geq n+1}\{0,1\}$. 
To any $f\in L^0(S^n,S_{2^n})$ associate the element $\tilde f\in [G]$,
defined for $x= (x_k)_{k\geq 1}\in S$, by
$\tilde f(x) = (f(x^n)(x_1,\ldots,x_n),x^n)$, where 
$x^n = (x_{n+1},x_{n+1},\ldots)\in S^n$. 
Let $L_n$ denote the subgroup of $[G]$ formed by 
$\{\tilde f\in [G]\colon f\colon S^n\to {\mathfrak{S}}_{2^n}\mbox{ is measurable}\}$.
Let $\tilde\mu$ be the probability measure on $S^n$, given for $A\subset S^n$,
$A$ measurable, by
\[\tilde\mu(A)=\sum_{\bar x\in S_n} \mu\left(\{ (\bar x,y)\colon y\in A
\}\right).\]
For $f\in L^0(S^n,{\mathfrak{S}}_{2^n})$ we have
\begin{eqnarray*}
\mu\left(\{x\in S\colon \tilde f(x)\neq x\} \right)&=&
\sum_{\bar x\in S_n}\mu\left(\{ (\bar x,y)\colon f(y)(\bar x) 
\neq \bar x\} \right) \\[2mm]
&\subset& \sum_{\bar x\in S_n}\mu\left(\{ (\bar x,y)\colon 
f(y)\neq {\mathrm{Id}}\}\right) \\[2mm]
&=& \tilde\mu\left(\{y\in S^n\colon f(y)\neq {\mathrm{Id}} \} \right).
\end{eqnarray*}
Therefore, if $L^0(S^n,{\mathfrak{S}}_{2^n})$ is endowed with the topology of convergence in measure
with respect to $\tilde\mu$ and $L_n\subset [G]$ with the uniform topology,
then the group isomorphism $L^0(S^n,{\mathfrak{S}}_{2^n})\ni f\mapsto \tilde f\in L_n$
is continuous. 

By Theorem \ref{gener}, $(L_n)_{n\geq 1}$ forms a sequence of extremely
amenable subgroups of $[G]$, with $[G_n]\subset L_n$,
$\forall n\geq 1$. As $([G_n])_{n\geq 1}$ forms an increasing sequence of
subgroups of $[G]$, whose union is uniformly dense, we conclude by Lemma
\ref{dense} that if $G$ and $X$ are as in Prop. \ref{hard},
then the full group $[G]$, endowed with the uniform
topology, is extremely amenable. 

Notice, however, that the actual conclusion
of our Proposition \ref{hard} is stronger, because there exist extremely
amenable groups that are not L\'evy groups, cf. e.g. Remark \ref{eanonl}.
\end{remark}

\subsection{} In this subsection, we show that if $[{\mathcal R}]$ 
endowed with the
uniform topology is amenable as a topological group, then $\mathcal R$ is an
amenable equivalence relation.

If $E$ is a separable Banach space, let $\Iso(E)$ denote the group of
isometric isomorphisms of $E$, endowed with the Borel structure induced by
the strong operator topology. As is well-known (see e.g.
\cite{zimmer-jfa}, Lemma 1.3), the map from
$\Iso(E)$ to $\Homeo(E_1^\ast)$, given by duality, is continuous, hence Borel,
if we consider the topology of uniform convergence on $\Homeo(E_1^\ast)$. 

Let $\alpha\colon {\mathcal{R}}\to\Iso(E)$
be a cocycle, that is, a measurable map such that
\[\alpha(x,y)\alpha(y,z)=\alpha(x,z)~~\mbox{ if }~~
(x,y), (y,z)\in{\mathcal{R}}.\]

Let $\alpha^\ast\colon {\mathcal{R}}\to \Homeo(E_1^\ast)$ denote the cocycle
given by
\[\alpha^\ast(x,y) =(\alpha(x,y)^{-1})^\ast =\alpha(y,x)^\ast,~~
(x,y)\in{\mathcal{R}}.\]

For $f\in L^1(S,E)$ and $\sigma\in [{\mathcal R}]$, 
let $T(\sigma)f$ denote the element
of $L^1(S,E)$ given by
\[T(\sigma)f(s) = r(s,\sigma)\alpha(s,s\sigma)f(s\sigma),~~
s\in S,\]
where $r(s,\sigma)$ denotes the Radon-Nikodym cocycle given for $s\in S$
and $\sigma\in [{\mathcal R}]$ by $d\mu(s\sigma) = r(s,\sigma)d\mu(s)$.

One checks easily that $T$ defines a representation of
$[{\mathcal R}]$ into $\Iso(L^1(S,E))$. 

\begin{lemma}
If $[{\mathcal R}]$ is endowed with the uniform topology, then 
$T$ is strongly continuous.
\end{lemma}

\begin{proof} If $\xi\in E$ and $f\in L^1(S,\C)$, let $\tilde f$ denote
the element of $L^1(S,E)$ given by
\[\tilde f(x) = f(x)\xi.\]
As the set of all such functions $\tilde f$ topologically spans $L^1(S,E)$,
it is enough to verify the continunity of every map of the form
\[[{\mathcal R}]\ni\sigma\mapsto T(\sigma)\tilde f\in L^1(S,E).\]
For $\sigma\in [{\mathcal R}]$ one has
\begin{eqnarray*}
\norm{T(\sigma)\tilde f - \tilde f} &=&
\int_S \norm{r(s,{\sigma})\alpha(s,s{\sigma})
\tilde f(s{\sigma}) - \tilde f(s)} d\mu(s)\\
&=&
\int_{\supp(\sigma)} \norm{r(s,\sigma) f(s{\sigma})\alpha(s,s{\sigma})\xi
- f(s)\xi}d\mu(s)\\
&\leq& \norm\xi \int_{\supp(\sigma)}\left(r(s,\sigma)\abs{f(s{\sigma})}+
\abs{f(s)}\right)d\mu(s)\\
&\leq& 2\norm\xi\norm f_\infty\mu(\supp(\sigma)) =
2\norm\xi\norm f_\infty d_\mu(\sigma,{\mathrm{Id}}),
\end{eqnarray*}
where $\supp(\sigma)=\{s\in S\colon s\sigma \neq s\}$. The statement follows.
\end{proof}

\subsection{}
Let ${\mathfrak{A}}=
\{A_s\colon s\in S\}$ be a Borel field of (non-empty) compact convex
subsets of $E^\ast$ such that for all $(x,y)\in{\mathcal{R}}$,
\[\alpha^\ast(x,y)(A_{y}) = A_x~~\mbox{ $\mu$-a.e.}\]
By \cite{zimmer-jfa}, Prop. 2.2, 
\[A=\{\lambda\in L^\infty(S,E^\ast)\colon \lambda(s)\in A_s ~~\mbox{ $\mu$-a.e.}
\}\]
is a weak$^\ast$ compact convex non-empty subset of $L^\infty(S,E^\ast)_1$.

Moreover, if $T^\ast$ denotes the dual representation of $T$ in 
$L^\infty(S,E^\ast)$, given for $\sigma\in [{\mathcal{R}}]$ and 
$\lambda\in L^\infty(S,E^\ast)$ by
\[T^\ast(\sigma)\lambda(s) =\alpha^\ast(s,s\sigma)f(s\sigma),
~~ s\in S,\]
then $A$ is invariant under the (affine) action $T^\ast$ of $[{\mathcal{R}}]$. 

We have already noted that this action is continuous if $A$ is equipped
with the weak$^\ast$ topology. Therefore, if the full 
group $[{\mathcal{R}}]$, equipped
with the uniform topology, is amenable as a topological group, then there
exists a $[{\mathcal{R}}]$-invariant section of the Borel field $\mathfrak A$. In
particular, this section is $\mathcal R$-invariant. Therefore, we have proved the
following proposition.

\begin{proposition} Let
$\mathcal R$ be a discrete measured equivalence relation, acting 
ergodically on a Lebesgue space $(S,\mu)$. If the full group $[{\mathcal{R}}]$,
equipped with the uniform topology, is amenable as a topological group,
then the equivalence relation $\mathcal R$ is amenable.
\label{ameme}
\end{proposition}

We summarize the results of this section in the following theorem.

\begin{theorem}
\label{criterion}
Let $\mathcal R$ be a discrete measured equivalence relation, acting 
ergodically on a Lebesgue space $(S,\mu)$. Then the following are equivalent.
\begin{enumerate}
\item \label{ame} $\mathcal R$ is an amenable equivalence relation.
\item \label{unea} $[{\mathcal{R}}]$ with the uniform topology is amenable.
\item\label{una}  $[{\mathcal{R}}]$ with the uniform topology is extremely
amenable.
\item\label{levy}  $[{\mathcal{R}}]$ with the uniform topology is a L\'evy group.
\end{enumerate}
\end{theorem}

\begin{proof} The implications (\ref{levy}) $\implies$
(\ref{una}) $\implies$ (\ref{unea})
are trivial, (\ref{unea}) $\implies$ (\ref{ame}) is Proposition \ref{ameme},
and (\ref{ame}) $\implies$ (\ref{levy}) is Proposition \ref{hard}.
\end{proof}

\begin{remark} Let $G$ be a discrete countable group acting freely and
ergodically on a Lebesgue space $(S,\mu)$.
Let $N[G]$ denote the normalizer of the full group.
Any element of $N[G]$ (resp. $[G]$) induces an automorphism
(resp. an inner automorphism) of the von Neumann factor $R$, normalizing
$M(L^\infty(S,\mu))$ (we keep the notation of Subs. \ref{fullgr}).
Then the $u$-topology (see Section \ref{agvna}) on $N[G]$ coincides with the
normal topology defined in Section 3.1 of \cite{danilenko} and is coarser
than the uniform topology. By Proposition \ref{hard}, if
$(S,\mu)$ is an amenable $G$-space, then $[G]$, endowed with the $u$-topology,
is a L\'evy group and in particular extremely amenable. 
\end{remark}

\section{Groups of non-singular transformations}
\subsection{Extreme amenability of $\Aut^\ast(X,\mu)$ with the weak topology}

In this subsection, we will prove the following result.

\begin{theorem} 
\label{th:aut*levy}
$\Aut^\ast(X,\mu)_w$ with the weak topology is a
L\'evy group and therefore extremely amenable. 
\end{theorem} 

If $(X, {\mathcal F}, \mu)$ is a non-atomic Lebesgue measure space, then it is
isomorphic to the circle $\s^1$, with the Lebesgue measure $\lambda$ and
therefore the topological groups
$\Aut^\ast(X,\mu)$ and $\Aut^\ast(\s^1,\lambda)$ are isomorphic. 
In \cite{katznelson}, Katznelson shows (Th. 4.1) that 
the set of all orientation preserving
$C^\infty$ diffeomorphisms of the circle $\s^1$ of type $III_1$ is dense in the 
group of all measure class preserving $C^\infty$ diffeomorphisms of $\s^1$.
In particular this insures the existence of ergodic type $III_1$
transformations of $(X,\mu)$.

Theorem \ref{th:aut*levy} will therefore 
follow from Proposition \ref{hard} and the 
following Proposition.

\begin{proposition} It $T$ is an ergodic measure class preserving
transformation of $(X,\mu)$ of type $III_1$, then the full group $[T]$ is 
weakly dense in $\Aut^\ast(X,\mu)$.
\label{p:III}
\end{proposition}



In order to prove Prop. \ref{p:III}, we need to introduce a few concepts.

Following Tulcea \cite{IT}, call a measurable mapping
$f\colon X\to Y$ between two finite measure spaces a
{\it Lebesgue mapping} if for every measurable $A\subseteq X$
\[\mu_Y(f(A))=\frac{\mu_Y(Y)}{\mu_X(X)}\mu_X(A).\]
If $X$ and $Y$ are measurable subsets of a standard Lebesgue probability space,
there is always a Lebesgue mapping from $X$ onto $Y$.

\begin{definition}
\label{permutation}
Say that a non-singular
transformation $g\in\Aut^\ast (X,\mu)$ is a {\it Lebesgue permutation} with
respect to a finite measurable partition $\mathcal P$ of $X$, if 
\begin{enumerate}
\item for each element $A\in{\mathcal P}$ the image
$g(A)$ also belongs to $\mathcal P$, and the restriction 
$g\vert_A$ is a Lebesgue mapping, and 
\item  if $A\in {\mathcal P}$ and $g(A)=A$, then
$g\vert_A=\mathrm{Id}_A$.
\end{enumerate}
\end{definition}

The following was deduced by Tulcea 
(\cite{IT}, Theorem 1) from a theorem, attributed by her to Linderholm.

\begin{theorem}
Lebesgue permutations are everywhere dense 
in the uniform topology on the group
$\Aut^\ast (X,\mu)$. Moreover, given a transformation $\tau\in\Aut^\ast (X,\mu)$,
an $\e>0$, and a 
finite measurable partition $\mathcal Q$ of $X$, there exists a finite
partition ${\mathcal Q}_1\prec{\mathcal Q}$ and a
Lebesgue permutation $\sigma$
with respect to ${\mathcal Q}_1$ such that $d_{unif}(\tau,\sigma)<\e$.
\label{tulcea}
\qed
\end{theorem}

\begin{proof}[Proof of Proposition \ref{p:III}]
Let $V$ be a weak neighbourhood of the identity of $\Aut^\ast(X,\mu)$,
determined by an $\e>0$ and a finite measurable partition
$\mathcal Q$ of $X$, so that
\[V=\{\beta\in\Aut^\ast (X,\mu)\colon \norm{\chi_A-\beta\chi_A}_1<\e
\mbox{ for each } A\in{\mathcal Q}\}.\]

By Theorem \ref{tulcea}, it is enough, given 
a finite partition ${\mathcal P}=\{P_1,P_2,\ldots,P_n\}\prec 
{\mathcal Q}$ of $(X,\mu)$ and a Lebesgue permutation $\tau$ 
of $(X,\mu)$ with respect to ${\mathcal P}$, to find a $\sigma\in [T]$
with $\sigma^{-1}\tau\in V$.

By Lemma 32 in \cite{HO}, for each $1\le i \leq n$, there exists 
$\sigma_i\in [T]$ such that
\begin{enumerate}
\item $\sigma_i(P_i) =\tau(P_i)$ and if $\tau(P_i)=P_i$, then
$\sigma_i=\mathrm{Id}$.
\item $\left\vert  \frac {d\mu\circ \sigma_i}{d\mu}(x)  
- \frac {\mu(\tau(P_i))}{\mu(P_i)} \right\vert < \varepsilon$
for $\mu$-a.e. $x \in P_i$. 
\end{enumerate}
Let $\sigma\in [T]$ be given by $\sigma\vert_{P_i}=\sigma_i$. 
As a consequence of (2), one has for $\mu$-a.e. $x\in X$,
\[ \left \vert\frac{d\mu\circ \sigma}{d\mu}(x) -\frac{d\mu\circ \tau}{d\mu}(x)
\right \vert<\e.\]
Moreover as ${\mathcal P}\prec {\mathcal Q}$, one has 
$\sigma^{-1}\tau(A)=A$
for each $A\in {\mathcal Q}$, and so $\sigma^{-1}\tau\in V$, as required.
\end{proof}

\begin{remark}
We do not know if the group $\Aut^\ast(X,\mu)_u$ with the uniform topology
is (extremely) amenable. 
\end{remark}

\subsection{Groups of isometries of $L^p$} 
In this subsection we will prove the following result.

\begin{theorem} Let $(X,\mu)$ be a standard Borel space with a non-atomic
measure and let $1\leq p<\infty$. 
The group of isometries $\Iso(L^p(X,\mu))$ with the strong operator
topology is a L\'evy group (and therefore extremely amenable).
\label{grisom}
\end{theorem}

In case $p=2$, the result by Gromov and Milman \cite{GrM} that 
the unitary group $U(\ell^2)$ with the strong operator topology 
is a L\'evy group has started the current problematics. For $p\neq 2$,
however, the proof is different. 

A description of
the isometries of $L^p(X,\mu)$ was obtained by S. Banach
(\cite{banach}, th\'eor\`eme 11.5.I, p. 178), but the proof was apparently 
never published. The first available proof (of a more general result
concerning not necessarily surjective isometries) belongs to 
Lamperti \cite{lamperti}, whose statement of
the main result (Theorem 3.1) was not quite correct, as noticed in
\cite{cfg}, cf. also \cite{grz}. 

\begin{theorem}[Banach] Let $1\leq p<\infty$, $p\neq 2$, let $(X,\mu)$ be a
finite measure space, and let
$T$ be a surjective isometry of $L^p(X,\mu)$. Then there is
an invertible measure class preserving transformation $\sigma$ of $X$ and a
measurable function $h$ with $\abs h =1$ such that
\begin{equation}
Tf = h\,\,^\sigma\! f.
\end{equation}
\qed
\label{lampcor}
\end{theorem}

%
%
%
%

Recall that $\Aut^\ast(X,\mu)_w$ can be identified with a (closed) 
topological subgroup of the group $\Iso(L^p(X,\mu))$,
equipped with the strong operator topology, through the left quasi-regular
representation:
\[^\sigma \!f(x) = f(\sigma^{-1}x)
\left(\frac{d(\mu\circ \sigma^{-1})}{d\mu}(x)\right)^{\frac 1p},
~~\sigma\in \Aut^\ast(X,\mu).\] 
Invertible (the same as onto) isometries of $L^p(X,\mu)$ correspond to the
invertible regular set isomorphisms, that is, invertible measure class 
preserving transformations of $(X,\mu)$, in which case 
$h\phi(f) = \,^{\phi^{-1}}\!\!f$.
We obtain:

\begin{corollary} Let $1\leq p<\infty$, $p\neq 2$, let $(X,\mu)$ be a
finite measure space. The group $\Iso(L^p(X,\mu))$ is the
semidirect product of the subgroup $\Aut^\ast(X,\mu)$ and the normal subgroup
$L^0(X,\mu;U(1))$ of all measurable maps from $X$ to the circle rotation group.
Moreover, the group $\Iso(L^p(X,\mu))$ with the strong operator topology is
the semidirect product of $\Aut^\ast(X,\mu)_w$ and the group
$L^0(X,\mu;U(1))$ equipped with the topology of convergence in measure.
\label{semidirect}
\end{corollary}

\begin{proof} According to Banach's theorem \ref{lampcor}, every element
$T\in \Iso(L^p(X,\mu))$ admits a (clearly unique) decomposition of the
form $Tf = h\cdot\,^\sigma\! f$, for $f\in L^p(X,\mu)$,
where $h\in L^0(X,\mu;U(1))$ and 
$\sigma\in \Aut^\ast(X,\mu)$. To establish the first claim, 
it is therefore enough to prove that 
$L^0(X,\mu;U(1))$ is normal in $\Iso(L^p(X,\mu))$ or, which is the same, 
invariant under inner automorphisms generated by the 
elements of $\Aut^\ast(X,\mu)$. Let $h\in L^0(X,\mu;U(1))$ and 
$\sigma\in \Aut^\ast(X,\mu)$. One has
\begin{eqnarray*}
(\sigma^{-1}h\sigma) (f)(x)  &\equiv & 
\,\,^{\sigma^{-1}}\!\!\left( h\cdot \,\,^{\sigma}\! f\right)(x) \\
&=& \left( h\cdot \,\,^{\sigma} \!f\right)(\sigma x)
\left(\frac{d(\mu\circ \sigma)}{d\mu}(x)\right)^{\frac 1p} \\
&=& h(\sigma x)f(\sigma^{-1}\sigma x)
\left(\frac{d(\mu\circ \sigma^{-1})}{d\mu}(\sigma x)\right)^{\frac 1p}
\left(\frac{d(\mu\circ \sigma)}{d\mu}(x)\right)^{\frac 1p}
\\ &=& \sigma^{-1}(h)(x)f(x),
\end{eqnarray*}
and therefore $\sigma^{-1}h\sigma = \sigma^{-1}(h)\in L^0(X,\mu;U(1))$.

We have already noted that 
$\Aut^\ast(X,\mu)_w$ is a closed topological subgroup of $\Iso(L^p(X,\mu))$.
The same is true of $L^0(X,\mu;U(1))$. Indeed, every $L^p$-metric,
$1\leq p<\infty$, induces the topology of convergence in measure on 
$L^0(X,\mu;U(1))$ (See e.g. Ex. 13.33 in \cite{HS}.)
Using this observation, it is easy to verify that,
first, for every $f\in L^p(X,\mu)$ the orbit map of the multiplication action
\[L^0(X,\mu;U(1))\ni h\mapsto h\cdot f\in L^p(X,\mu)\]
is continuous, while for $f(x)\equiv 1$ it is a homeomorphic embedding.

As a corollary (which can be also checked directly),
the action of $\Aut^\ast(X,\mu)_w$ on $L^0(X,\mu;U(1))$, defined by
$(\sigma, h)\mapsto \sigma^{-1}(h)$, is continuous as
a map $\Aut^\ast(X,\mu)_w\times L^0(X,\mu;U(1))\to L^0(X,\mu;U(1))$.
Therefore, the semi-direct product 
$\Aut^\ast(X,\mu)_w \ltimes L^0(X,\mu;U(1))$ is a (Polish) topological
group. 

The algebraic automorphism given by the multiplication map,
\[\Aut^\ast(X,\mu)_w \ltimes L^0(X,\mu;U(1))\to \Iso(L^p(X,\mu)),\]
is continuous, and since all groups involved are Polish, it is a topological
isomorphism by the standard Open Mapping Theorem for Polish groups, 
cf. \cite{husain}, Corollary 3, p. 98 and Theorem 3, p. 90.
\end{proof} 

\begin{proof}[Proof of Theorem \ref{grisom}]
As in the Subsection \ref{eafg}, we will identify the measure space
$(X,\mu)$ with the compact abelian group $X=\prod \{0,1\}$ equipped with an
ergodic non-atomic quasi-invariant measure $\mu$.
Let $L_n$ denote the subgroup of $L^0(X,\mu;U(1))$ consisting
of functions constant on cylinder sets determined by the first $n$ coordinates.
It is clear that the
group $L_n$ is invariant under all inner automorphisms generated by
elements of the permutation group $\Bij(n)$ of the finite set $S_n$.

The normalized Haar measures on $\Bij(n)$ concentrate in $\Aut^\ast(X,\mu)_w$  by 
Corollary \ref{largestatom}. Since concentration depends on the uniform
structure only, one can assume without loss in generality that the groups
$L_n$ are equipped with the $L^1(\mu)$-metric induced from $L^0(X,\mu;U(1))$
(because it induces the topology of convergence in measure on the former
group and is translation invariant, therefore induces the additive uniformity).
By Corollary \ref{largestatom1}, the normalized Haar measures on $L_n$,
$n\in\N$, concentrate in $L^0(X,\mu;U(1))$.

The semidirect product groups $\Bij(n)\ltimes L_n$ are 
compact, and the normalized Haar measures on them are products of Haar
measures, respectively, on $\Bij(n)$ and on $L_n$. Therefore, those measures
concentrate in $\Iso(L^p(X,\mu))$ by Lemma \ref{subadd}. 

Finally, the union of subgroups $\Bij(n)\ltimes L_n$ is everywhere dense in
$\Iso(L^p(X,\mu))$ by Corollary \ref{semidirect}, because the union
of $\Bij(n)$ is uniformly
dense in its own full group by Lemma \ref{approximation},
and therefore weakly dense in $\Aut^\ast(X,\mu)$ by Proposition \ref{p:III},
while the groups $L_n$ of simple functions are dense in $L^0(X,\mu;U(1))$.
\end{proof}

Theorem \ref{grisom} was conjectured by one of the present authors
in \cite{P02a}.

\section{Concentration to a nontrivial space}

M. Gromov \cite{Gr}, section 3$\frac 12$.45 (p. 200) has defined a
metric on the space of equivalence classes of Polish 
$mm$-spaces in such a way
that a sequence $X_n=(X_n,d_n,\mu_n)$ of $mm$-spaces forms a 
L\'evy family if and only if
it converges to the trivial $mm$-space $\{\ast\}$.

This approach allows one to talk of {\it concentration to a nontrivial 
$mm$-space.}
According to Gromov, this type of concentration commonly occurs
in statistical physics. At the same time, there are very few known
non-trivial examples
of this kind in the context of transformation
groups. The aim of this Section is to give some natural 
examples of concentration of this type.

We will recall Gromov's construction.
Let $\mu^{(1)}$ denote the
Lebesgue measure on the unit interval $\I=[0,1]$.
On the space $L^0(0,1)$ of all (equivalence classes) of
measurable real-valued functions define the metric
$\operatorname{me}_1$, generating the topology of convergence in measure,
as follows: 
$\operatorname{me}_1(h_1,h_2)$ is the infimum of all $\lambda>0$ with
the property 
\[
\mu^{(1)}\{\vert h_1(x)-h_2(x)\vert>\lambda\}<\lambda.
\]

Let $X=(X,d_X,\mu_X)$ and
$Y=(Y,d_Y,\mu_Y)$ be two Polish $mm$-spaces. 
There exist measurable maps
$f\colon\I\to X$, $g\colon\I\to Y$ such that $\mu_X=f_\ast\,\mu^{(1)}$
and $\mu_Y=g_\ast\,\mu^{(1)}$. Denote by $L_f$
the set of all functions of the form $h= h_1\circ f$, where
$h_1\colon X\to \R$ is 1-Lipschitz, having the additional property $h(0)=0$.
Similarly, define the set $L_g$. Both $L_f$ and $L_g$ are compact
subsets of $L^0(0,1)$. 

Let
$\HL(X,Y)$ be the infimum of Hausdorff distances between $L_f$
and $L_g$ (with regard to the metric $\operatorname{me}_1$), 
taken over all parametrizations $f$ and $g$ as above.
One can prove that $\HL$ is a metric on the space of (isomorphism classes
of) all Polish $mm$-spaces. (The difficult part of the proof is verifying
the first axiom of a metric.)

A sequence of $mm$-spaces
$X_n=(X_n,d_n,\mu_n)$ forms a L\'evy family if and only if
it converges to the trivial (one-point) $mm$-space in the metric $\HL$:
\[
X_n\overset{\HL}{\longrightarrow} \{\ast\}.
\]
(This is a reformulation of Lemma \ref{nearlyconstant}.)

 \begin{theorem}[Gromov] Let $Y=(Y,d_Y,\nu)$ and $Z_n=(Z_n,d_n,\mu_n)$, 
 $n\in\N$ be $mm$-spaces, where $(Z_n)$ forms a L\'evy family.
 Let $d_Y\oplus^{(2)} d_n$ denote the $\ell_2$-type sum of the metrics $d_Y$ and
$d_n$ on the product space 
 $Y\times Z_n$.
 Then the family of $mm$-spaces $(Y\times Z_n, d_Y\oplus^{(2)} d_n, \nu\otimes\mu_n)$
 concentrates to the $mm$-space $(Y,d_Y,\nu)$:
 \[(Y\times Z_n, d_Y\oplus^{(2)} d_n, \nu\otimes\mu_n)
 \overset{\HL}{\longrightarrow}(Y,d_Y,\nu).\]
 \label{gromov}
 \qed
 \end{theorem}

We are going to show a situation where 
Gromov's theorem applies naturally.

\begin{lemma} 
Let $G$ be a metrizable topological group which is, as a
topological group, a semidirect product
of a compact subgroup $H$ and a closed normal subgroup $N$:
\[G = H\ltimes N.\]
Assume that $N$ is a L\'evy group, that is, there is a 
directed family of compact subgroups $(K_\alpha)$ of $N$, having
an everywhere dense union and such that $(K_\alpha, \mu_\alpha)$ form 
a L\'evy family with regard to the right
uniform structure on $G$. 
Let $\nu$ denote
the normalized Haar measure on $H$. 
Then there exist a right-invariant compatible metric $\rho$ on $G$ 
such that
the family of $mm$-spaces $H\cdot K_\alpha$, equipped with the metric $\rho$
and the convolution of normalized Haar measures, concentrates 
to the $mm$-space $(H,\rho\vert_H,\nu)$:
\[(H\cdot K_\alpha,\rho,\nu\ast\mu_\alpha)\overset{\HL}{\longrightarrow} 
(H,\rho\vert_H,\nu).\]
If in addition each of the subgroups $(K_\alpha)$ is invariant under the action
of $H$ by conjugations, then every $H\cdot K_\alpha$ is a compact subgroup of
$G$.
\label{7.3}
\end{lemma}

\begin{proof} Denote by $\tau$ the continuous
action of $H$ on $N$ by topological group isomorphisms, so that
$G = H\ltimes_\tau N$. Fix any compatible right-invariant metric, $d$,
on $N$, and define a new metric, $\varsigma$, as follows: for
$x,y\in N$, set
\[\varsigma_N(x,y) =\int_H d(\tau_hx,\tau_hy)~d\nu(h).\]
This $\varsigma_N$ is also a right-invariant compatible metric on $N$, which
in addition is invariant under the action $\tau$.

Now choose a bi-invariant compatible metric, $\varsigma_H$,
on the compact group $H$. Define a metric $\rho$ on $G$ as the $\ell_2$-type 
sum of $\varsigma_H$ and $\varsigma_N$:
\[\rho((h,x),(h',x')) = \sqrt{\varsigma_H(h,h')^2+\varsigma_N(x,x')^2}.\]

It is easily verified that $\rho$ is a right-invariant compatible metric on
$G$. Besides, $\rho$ is invariant under left multiplication by elements of
$H$:
\begin{eqnarray*}
\rho((h,e)\cdot(g,x),(h,e)\cdot(g',x')) &=& 
\rho((hg,\tau_hx),(hg',\tau_hx')) \\
&=& \left(\varsigma_H(hg,hg')^2+\varsigma_N(\tau_hx,\tau_hx')^2
\right)^{\frac 12} \\
&=& \left(\varsigma_H(g,g')^2+\varsigma_N(x,x')^2 
\right)^{\frac 12} \\
&=& \rho((g,x),(g',x')).
\end{eqnarray*}

The restriction of $\rho$ to each $H\cdot K_\alpha$ is just the
$\ell_2$-type sum of $\varsigma_H$ with the restriction of $\varsigma_N$ to
$K_\alpha$.

Since the multiplication mapping
$H\times K_\alpha\to H\cdot K_\alpha$ is one-to-one (by assumption) and 
also of course continuous, it is a homeomorphic embedding of 
$H\times K_\alpha$
into our topological group $G$. In particular, the convolution of measures
$\nu\ast\mu_\alpha$ is in fact the product measure, $\nu\otimes\mu_\alpha$.
Now Gromov's theorem \ref{gromov} applies. 

The last statement of the Lemma is obvious.
\end{proof} 

\begin{remark} As a direct consequence from \cite{P02a}, Lemma 2.3,
under assumptions of Lemma \ref{7.3}, $H$ forms the universal
minimal flow of the topological group $G$. 
It would be interesting to find
more general situations where concentration of subobjects of $G$ to a 
nontrivial compact
space $X$ would mean that $X$ is a universal minimal flow of $G$.
\label{unimin}
\end{remark}

Using the description of the automorphism groups of injective factors of
type $III$, we now present concrete examples of concentration to a
non-trivial space. In particular, in the last Example \ref{sbgrps} we get
a sequence of compact subgroups concentrating to a non-trivial space with
regard to one group topology, and to a trivial space with regard to
another group topology.

\begin{example}
\label{rlambda}
For $0<\lambda<1$, let us consider $G=\oplus (\Z/2)$ acting freely and
ergodically on $(X,\mu_\lambda) =
\prod\left(\{0,1\},\left\{\frac 1{1+\lambda},\frac\lambda{1+\lambda}\right\}
\right)$, and let $R_\lambda$ denote the von Neumann factor of type
$III_\lambda$, obtained by the group measure space construction. 
Recall that $R_\lambda$ is approximately finite dimensional and by 
\cite{connes-injective} is up to isomorphism the unique injective factor of
type $III_\lambda$. To any automorphism $\alpha$ of $R_\lambda$ is
associated its modulus, $\mod\alpha$, which belongs to
$\R/(\log\lambda)\Z$ (see for example \cite{Tak-3}, Chap. XII). 
If $\Aut(R_\lambda)$ is endowed with the $u$-topology, then the map
$\mod\colon \Aut(R_\lambda)\to\R/(\log\lambda)\Z$ is a continuous
homomorpism. Moreover (see \cite{Tak-3}, Thm. XVIII.1.11) there is a natural
short exact sequence 
\[1\to \overline{\Inn(R_\lambda)}\to \Aut(R_\lambda) 
\overset{\mathrm{mod}}{\to} \R/(\log\lambda)\Z\to 1,\]
and a continuous action $\tilde\theta\colon \R/(\log\lambda)\Z\to
\Aut(R_\lambda)$ with 
\[\mod(\tilde\theta(s))=s,~~ s\in \R/(\log\lambda)\Z.\]
By Proposition \ref{pcn}, $\Inn(R_\lambda)$ and consequently
$\overline{\Inn(R_\lambda)}$ are L\'evy groups.
Hence by Lemma \ref{7.3} the group $\Aut(R_\lambda)$ contains a chain of
compact $mm$-spaces with everywhere dense union which
concentrates to $\R/(\log\lambda)\Z$.
\end{example}

\begin{example}
Let $G$ be a countable discrete group acting freely and ergodically on
the Lebesgue space $(X,\mu)$. By Proposition \ref{hard}, if $(X,\mu)$ is an
amenable $G$-space, then the full group $[G]$, endowed with the $u$-topology,
is extremely amenable.

If the dynamical system $(X,\mu,G)$ is of type $III_\lambda$,
$0<\lambda<1$, then there is a split exact sequence 
\[1\to \overline{[G]}^{u{\mathrm{-top}}}\to N[G]\overset{\mathrm mod}{\to}
 \R/(\log\lambda)\Z\to 1\]
(See \cite{danilenko}, 3.3.)
As $\overline{[G]}^{u{\mathrm{-top}}}$ is extremely amenable and moreover a 
L\'evy group, $N[G]$ contains a chain of
compact $mm$-spaces with everywhere dense union which
concentrates to $\R/(\log\lambda)\Z$.
\end{example}

\begin{example}
\label{sbgrps}
Let $R_\lambda$ be the injective factor of type $III_\lambda$ (see Example
\ref{rlambda}), realized as an infinite tensor product of factors of type
$I_2$, i.e. $R_\lambda=\otimes (M_2(\C),\varphi_\lambda)$, where  
$\varphi_\lambda = {\mathrm{Tr}}(h_\lambda,\cdot)$ and $h_\lambda$ is the 
density matrix $\left(\begin{matrix} \frac 1 {1+\lambda} &  \\ & \frac \lambda
 {1+\lambda}   \end{matrix}\right)$. If $\varphi$ denotes the state
 $\otimes \varphi_\lambda$, then the modular group of automorphisms
 $\sigma^{\varphi}$ is given by
 \[\sigma_t^{\varphi} =\otimes {\mathrm{Ad}}\, 
 \left(\begin{matrix} 1 &  \\ & \lambda^{i t}   \end{matrix}\right),
 ~~t\in\R.
\]
For $m\geq 1$, let $N_m$ be the type $I_{2^m}$ subfactor of $M$ generated by
\[\{x_1\otimes x_2 \otimes \ldots \otimes x_m \otimes 1 \otimes 1 \otimes
\ldots \colon x_k\in M_2(\C)\}.\]
Let $K_m$ denote the compact subgroup of $\Inn(R_\alpha)$ given by
$\{{\mathrm{Ad}}\,(u)\colon u\in U(N_m)\}$. 
Let $\Cnt(R_\lambda)\triangleleft\Aut(R_\lambda)$ be the normal subgroup of
centrally trivial automorphisms of $R_\lambda$. 
Let $T = \frac{- 2\pi}{\log\lambda}$.
By \cite{Tak-3}, Thm. XVIII.1.11, 
\begin{equation}
\Cnt(R_\lambda) =\Inn(R_\lambda)\rtimes_{\sigma^\varphi} \R/T\Z.
\label{semi}
\end{equation}
By Propositions
\ref{increasing} and \ref{pcn}, $\{K_m\}_{m\geq 1}$ forms a L\'evy family
of subgroups of $\Inn(R_\lambda)$. 
As for $t\in \R$, $\sigma^{\varphi}_t$ leaves $N_m$ globally
fixed, 
\[K_m\cdot \R/T\Z =
\{{\mathrm{Ad}}\,u\circ \sigma^{\varphi}_s\colon u\in U(N_m),~~s\in\R\}\]
is a subgroup of $\Cnt(R_\lambda)$. 

We will now consider two group topologies on $\Cnt(R_\lambda)$, with regard
to which the family of compact subgroups $(K_m\cdot \R/T\Z)_{m=1}^\infty$
exhibit different kinds of concentration behaviour.

The first one is the topology of semidirect product using the decomposition
in Eq. (\ref{semi}). Consider $\Inn(R_\lambda)$ as the
topological factor-group of the group $U(R_\lambda)$, equipped with the
strong operator topology (or the $s(R_\lambda,(R_\lambda)_\ast)$-topology.
As a factor-group of a Polish group by a compact subgroup (the central
subgroup isomorphic to the circle rotation group),
$\Inn(R_\lambda)$ is a Polish group itself.
It is easy to check that the action $\sigma^\varphi$ of $\R/T\Z$ on 
$\Inn(R_\lambda)$ is continuous. Therefore, $\Cnt(R_\lambda)$ is a
Polish group, and it contains a chain of
compact subgroups $K_m\cdot \R/T\Z$, $m\in\N$
with everywhere dense union, which family of subgroups
concentrate in Gromov's sense to the circle rotation group 
$\R/T\Z$ by force of Lemma \ref{7.3}. 

Notice also that the compact space $\R/T\Z$ forms the universal minimal flow of
the topological group $\Cnt(R_\lambda)$, cf. Remark \ref{unimin}.

The second topology on the group $\Cnt(R_\lambda)$
is the familiar $u$-topology. In this case,
$\Inn(R_\lambda)$ is everywhere dense in $\Cnt(R_\lambda)$, and in particular the 
group $\Cnt(R_\lambda)$ is extremely amenable. Let $d$ be a 
compatible metric on $\Cnt(R_\lambda)$, invariant on one side (e.g. on the 
left). 
For every $\e>0$, if $m$ is sufficiently large, the $\e$-neighbourhood
of $K_m$ contains the subgroup $\R/T\Z$ and therefore all of $\Cnt(R_\lambda)$.
Now one can easily see that
the compact subgroups $K_m\cdot \R/T\Z$, $m\in\N$ form a L\'evy family and
thus, in contrast to the previous situation, concentrate to a one-point space.
\end{example}

\end{document}